\documentclass{amsart}
\usepackage{amsmath,amssymb,graphicx}

\newtheorem{thm}{Theorem}
\newtheorem{lem}[thm]{Lemma}
\newtheorem{cor}[thm]{Corollary}

\theoremstyle{definition}

\newtheorem{rmk}[thm]{Remark}

\newcommand{\CPb}{\overline{\mathbb{CP}}{}^{2}}
\newcommand{\CP}{{\mathbb{CP}}{}^{2}}
\newcommand{\R}{\mathbb{R}}
\newcommand{\Z}{\mathbb{Z}}
\newcommand{\scparallel}{{\scriptscriptstyle \parallel}}

\title[Exotic small 4-manifolds with odd signatures]
{Exotic smooth structures on small\\ 
4-manifolds with odd signatures} 

\begin{document}

\author{Anar Akhmedov}
\address{School of Mathematics, 
University of Minnesota, 
Minneapolis, MN, 55455, USA}
\email{akhmedov@math.umn.edu}

\author{B. Doug Park}
\address{Department of Pure Mathematics, 
University of Waterloo, 
Waterloo, ON, N2L 3G1, Canada}
\email{bdpark@math.uwaterloo.ca}

\date{May 21, 2007.  Revised on July 1, 2009}


\begin{abstract}
Let $M$\/ be $\CP\#2\CPb$, $3\CP\#4\CPb$ or $(2n-1)\CP\#2n\CPb$ for any integer $n\geq 3$.
We construct an irreducible symplectic 4-manifold homeomorphic to $M$ and also an infinite family of pairwise non-diffeomorphic irreducible non-symplectic 4-manifolds homeomorphic to $M$.   
We also construct such exotic smooth structures when $M$\/ is $\CP\#4\CPb$ or $3\CP\# k \CPb$ for $k=6,8,10$.  
\end{abstract}

\maketitle

\section{Introduction}

This paper is a belated sequel to \cite{AP}.  A concise history of constructing exotic smooth structures on simply-connected $4$-manifolds with small Euler characteristics is found in the introduction of \cite{AP}.  Given two $4$-manifolds, $X$\/ and $Y$, we denote their connected sum by $X\# Y$.  For a positive integer $m\geq 2$, 
the connected sum of $m$\/ copies of $X$\/ will be denoted by $mX$\/ for short. 
Let $\CP$ denote the complex projective plane and let $\CPb$ denote the underlying 
smooth $4$-manifold $\CP$ equipped with the opposite orientation.    
Our main result is the following.

\begin{thm}\label{thm:main}
Let $M$ be one of the following\/ $4$-manifolds.  
\begin{itemize}
\item[(i)] $\CP\# m\CPb$ for $m=2,4$,

\vspace{2pt}
\item[(ii)] $3\CP\# k\CPb$ for $k=4,6,8,10$,  

\vspace{2pt}
\item[(iii)] $(2n-1)\CP\#2n\CPb$ for any integer $n\geq 3$.
\end{itemize}
Then there exist an irreducible symplectic\/ $4$-manifold and an infinite family of pairwise non-diffeomorphic irreducible non-symplectic\/ $4$-manifolds, all of which are homeomorphic to $M$.  
\end{thm}

Currently $\CP\#2\CPb$ has the smallest Euler characteristic amongst all simply-connected topological 4-manifolds that are known to possess more than one smooth structure.   Recall that exotic irreducible smooth structures on $(2n-1)\CP\#2n\CPb$ for $n\geq 49$ were already constructed in \cite{ABBKP}.  It is an intriguing open problem whether the exotic $4$-manifolds we construct in this paper (Section~\ref{sec:exotic 3-4}) in that range are in fact diffeomorphic to the corresponding $4$-manifolds in \cite{ABBKP}.  

Combined with the results in \cite{ABBKP}, Theorem~\ref{thm:main} allows us to conclude that for any pair of positive integers $(m,n)$ with $m$\/ odd and $m< n \leq 5m+4$, there is an irreducible symplectic 4-manifold homeomorphic to $m\CP\#n\CPb$ and moreover, there is an infinite family of pairwise non-diffeomorphic irreducible non-symplectic $4$-manifolds homeomorphic to such $m\CP\#n\CPb$.  By blowing up repeatedly, we also obtain such exotic smooth structures, albeit reducible, on $m\CP\#n\CPb$ with $m$\/ odd and $n>5m+4$.  
In terms of the geography problem, we conclude that there exist a simply-connected  irreducible symplectic $4$-manifold and infinitely many simply-connected pairwise non-diffeomorphic irreducible non-symplectic $4$-manifolds 
that realize the following coordinates:
\begin{itemize}
\item[] \hspace{9.4pt}$(e,\sigma)$ when $2e+3\sigma\ge 0$, $e+\sigma\equiv 0\pmod{4}$, and $\sigma\leq -1$, 

\vspace{3pt}
\item[]$(\chi_h,c_1^2)$ when $0 \leq c_1^2 \leq 8\chi_h -1$.

\vspace{2pt}
\end{itemize}
Here, $e$\/ and $\sigma$\/ denote the Euler characteristic and the signature respectively, while $\chi_h=(e+\sigma)/4$ and $c_1^2= 2e+3\sigma$.   The results in this paper fill in 52 new points in the geography plane that were left open in \cite{ABBKP}.  Combining Theorem~\ref{thm:main} above with Theorem~23 in \cite{ABBKP}, we can also deduce the following.  

\begin{thm}\label{thm:wedge}
Let\/ $X$ be a closed symplectic\/ $4$-manifold and suppose that\/ $X$ contains a symplectic torus\/ $T$ of self-intersection\/ $0$ such that the inclusion induced homomorphism\/ $\pi_1(T) \rightarrow \pi_1(X)$ is trivial.  Then for any pair\/ $(\chi, c)$ of nonnegative integers satisfying
\begin{equation}\label{eq: chi-c inequality}
0\leq c \leq 8\chi -1,
\end{equation}
there exists a symplectic\/ $4$-manifold\/ $Y$ with\/ $\pi_1(Y)=\pi_1(X)$, 
\begin{equation}\label{eq: chi_h and c_1^2}
\chi_h(Y)=\chi_h(X)+\chi \text{ \ and\/ \ } 
c_1^2(Y)=c_1^2(X)+c .
\end{equation}
Moreover, $Y$ has an odd indefinite intersection form, and if\/ $X$ is minimal then $Y$ is minimal as well.  
\end{thm}

Our paper is organized as follows.
In Sections~\ref{sec:(2n-3)(S^2 x S^2)}--\ref{sec: Luttinger surgery}, we collect building blocks that are needed in our construction of exotic $4$-manifolds.  
In Sections~\ref{sec: base point}--\ref{sec: pi_1}, we calculate the fundamental groups of some of our building blocks.  
In Section~\ref{sec:exotic 1-2}, we construct exotic smooth structures on $\CP\#2\CPb$.  
In Section~\ref{sec:exotic 3-4}, we construct exotic smooth structures on $(2n-1)\CP\# 2n\CPb$ for $n\geq 2$.  In Section~\ref{sec:exotic 1-4}, we construct exotic smooth structures on $\CP\#4\CPb$.   In Section~\ref{sec:exotic 3-k}, we construct exotic smooth structures on $3\CP\#k\CPb$ for $k=6,8,10$.  Finally in Section~\ref{sec:proof of wedge thm}, we present a proof of Theorem~\ref{thm:wedge}.

\section{Construction of cohomology $(2n-3)(S^2\times S^2)$}
\label{sec:(2n-3)(S^2 x S^2)}

\begin{figure}[ht]
\begin{center}
\includegraphics[scale=.49]{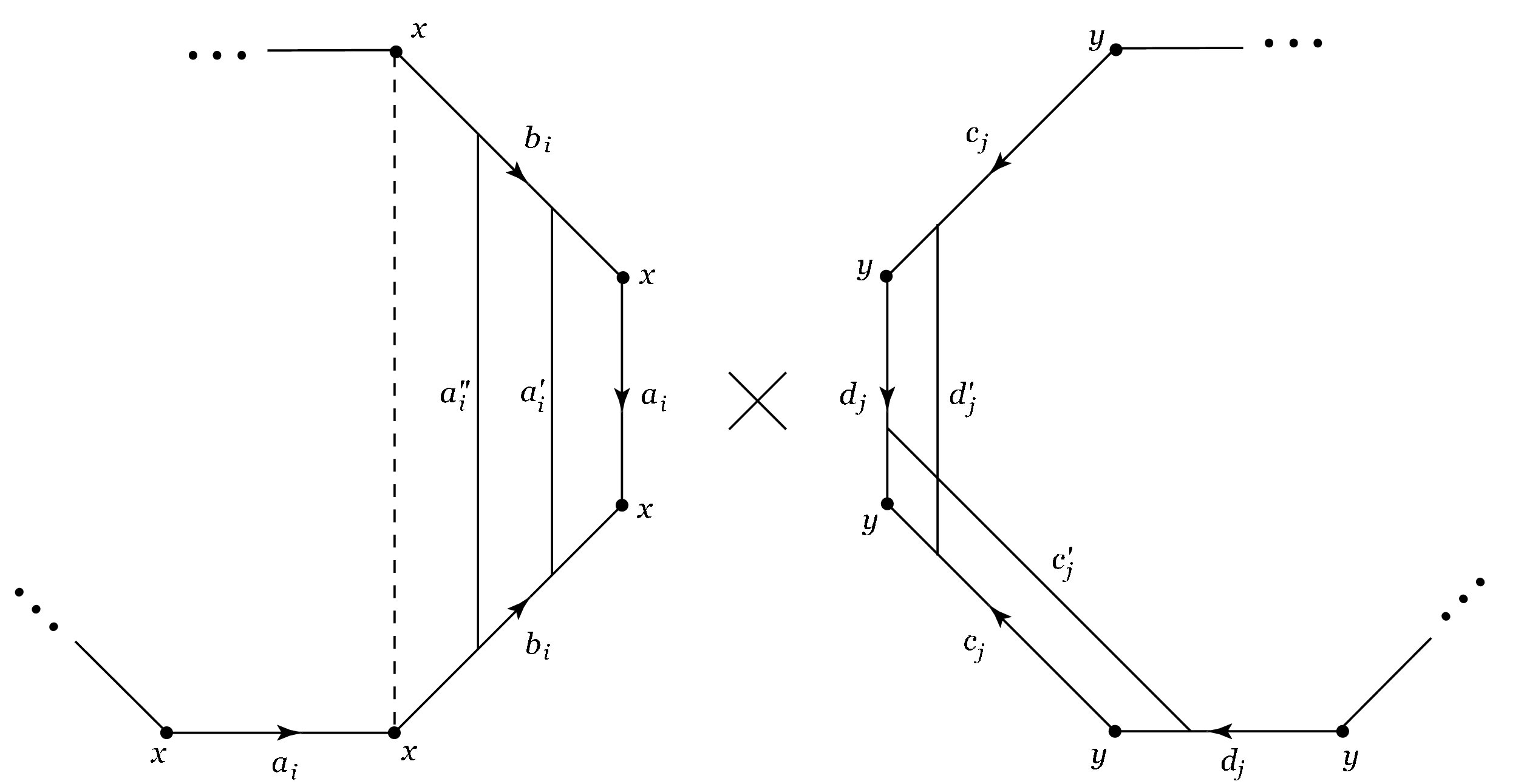}
\caption{Lagrangian tori $a_i'\times c_j'$ and $a_i''\times d_j'$}
\label{fig:lagrangian-pair}
\end{center}
\end{figure}

For each integer $n\geq2$, we construct a
family of irreducible pairwise non-diffeomorphic 4-manifolds
$\{Y_n(m)\mid m=1,2,3,\dots\}$ that have the same integer cohomology ring
as $(2n-3)(S^2\times S^2)$.  In what follows, we will use the notation of \cite{FPS} wherein the $n=2$ case is worked out in full detail.  $Y_n(m)$ are gotten by performing
$2n+3$ Luttinger surgeries (cf.\ \cite{ADK, luttinger}) and a single $m$\/ torus surgery on $\Sigma_2\times \Sigma_n$.  Here, $\Sigma_g$ denotes a closed Riemann surface of genus $g$.  
These $2n+4$ surgeries comprise of the following 8 torus surgeries in \cite{FPS}
\begin{eqnarray}\label{first 8 Luttinger surgeries}
&&(a_1' \times c_1', a_1', -1), \ \ (b_1' \times c_1'', b_1', -1), \ \
(a_2' \times c_2', a_2', -1), \ \ (b_2' \times c_2'', b_2', -1),\\ \nonumber
&&(a_2' \times c_1', c_1', +1), \ \ (a_2'' \times d_1', d_1', +1),\ \
(a_1' \times c_2', c_2', +1), \ \ (a_1'' \times d_2', d_2', +m),
\end{eqnarray}
together with the following $2(n-2)$ new Luttinger surgeries
\begin{gather*}
(b_1'\times c_3', c_3',  -1), \ \ 
(b_2'\times d_3', d_3', -1), \  \dots  ,\ 
(b_1'\times c_n', c_n',  -1), \ \
(b_2'\times d_n', d_n', -1).
\end{gather*}
Here, $a_i,b_i$ ($i=1,2$) and $c_j,d_j$ ($j=1,\dots,n$) are standard loops that generate $\pi_1(\Sigma_2)$ and $\pi_1(\Sigma_n)$, respectively.  The prime and double prime notations are explained in \cite{FPS}.  
Figure~\ref{fig:lagrangian-pair} depicts typical Lagrangian tori along which we perform surgeries.  

Recall that $Y_n(m)$ is symplectic only when $m=1$.  Also note that the Euler characteristic of $Y_n(m)$ is $4n-4$ and its signature is $0$.  $\pi_1(Y_n(m))$ is generated by $a_i,b_i,c_j,d_j$ ($i=1,2$ and $j=1,\dots,n$) and the following relations hold in $\pi_1(Y_n(m))$:  
\begin{gather}\label{Luttinger relations}
[b_1^{-1},d_1^{-1}]=a_1,\ \  [a_1^{-1},d_1]=b_1,\ \  [b_2^{-1},d_2^{-1}]=a_2,\ \  [a_2^{-1},d_2]=b_2,\\ \nonumber
[d_1^{-1},b_2^{-1}]=c_1,\ \ [c_1^{-1},b_2]=d_1,\ \ [d^{-1}_2,b^{-1}_1]=c_2,\ \ [c_2^{-1},b_1]^m=d_2,\\ \nonumber
 [a_1,c_1]=1, \ \ [a_1,c_2]=1,\ \  [a_1,d_2]=1,\ \ [b_1,c_1]=1,\\ \nonumber
[a_2,c_1]=1, \ \ [a_2,c_2]=1,\ \  [a_2,d_1]=1,\ \ [b_2,c_2]=1,\\ \nonumber
[a_1,b_1][a_2,b_2]=1,\ \ \prod_{j=1}^n[c_j,d_j]=1,\\ \nonumber
[a_1^{-1},d_3^{-1}]=c_3, \ \ [a_2^{-1} ,c_3^{-1}] =d_3, \  \dots, \ 
[a_1^{-1},d_n^{-1}]=c_n, \ \ [a_2^{-1} ,c_n^{-1}] =d_n,\\ \nonumber
[b_1,c_3]=1,\ \  [b_2,d_3]=1,\ \dots, \
[b_1,c_n]=1,\ \ [b_2,d_n]=1.
\end{gather}

The surfaces $\Sigma_2\times\{{\rm pt}\}$ and $\{{\rm pt}\}\times
\Sigma_n$ in $\Sigma_2\times\Sigma_n$ descend to surfaces in
$Y_n(m)$.  Let us call their images $\Sigma_2$ and $\Sigma_n$ for
short.  We have $[\Sigma_2]^2=[\Sigma_n]^2=0$ and
$[\Sigma_2]\cdot[\Sigma_n]=1$. They are symplectic surfaces in
$Y_n(1)$.  Let $\mu(\Sigma_2)$ and $\mu(\Sigma_n)$ denote the meridians
of these surfaces in $Y_n(m)$.

\section{Construction of a genus $2$ symplectic surface in $T^4\#\CPb$}
\label{sec:braided surface}

Recall that ${\rm Br}_2(D^2)\cong \Z$ (cf.\ \cite{birman, fadell-vanbuskirk}), where $D^2$ denotes a $2$-dimensional disk.  Let $\beta$\/ be a connected $2$-string braid that generates ${\rm Br}_2(D^2)$.   An embedding $D^2\hookrightarrow T^2$ allows us to view $\beta$\/ as an element of ${\rm Br}_2(T^2)$.  
Such $\beta$\/ gives rise to a smooth simple loop in the product $D^2 \times S^1 \subset T^2 \times S^1$, which we will also denote by $\beta$.  Note that the homology class of the closed curve $\beta$\/ is $2[ \{{\rm pt}\} \times S^1] \in H_1(T^2 \times S^1;\Z)$.  Thus we get an embedded torus $T_{\beta} = \beta \times S^1 \subset (D^2 \times S^1)\times S^1 \subset (T^2 \times S^1) \times S^1 = T^4$.  It is easy to see that $T_{\beta}$ is a symplectic submanifold of $T^4$ endowed with a product symplectic structure (cf.\ \cite{FS: symplectic surfaces}).  

Let $\alpha_i$, $i=1,\dots,4$, denote the standard circle factors that generate $\pi_1(T^4)\cong H_1(T^4;\Z)\cong\Z^4$.   The tori $\alpha_1\times\alpha_2$ and $\alpha_3\times \alpha_4$ are symplectic submanifolds with respect to our chosen product symplectic structure on $T^4=T^2\times T^2$.  Note that $[T_{\beta}]=2[\alpha_3\times\alpha_4]$ in $H_2(T^4;\Z)$.  Now $T_{\beta}$ intersects the symplectic torus $\alpha_1\times\alpha_2$ at two points.  Symplectically resolve one of the two intersection points to obtain a genus 2 surface with one positive double point.  Next, symplectically blow up at the double point to obtain a smooth genus 2 symplectic surface $\bar{\Sigma}_2$ in $T^4\#\CPb$.  The homology class of $\bar{\Sigma}_2$ in $T^4\#\CPb$ is given by the sum $[\alpha_1\times\alpha_2]+2[\alpha_3\times\alpha_4]-2[E]$, where $E$\/ is the exceptional sphere of the blow-up.  It is easy to see that the self-intersection of $\bar{\Sigma}_2$ is zero.  Let $\bar{a}_1,\bar{b}_1,!
 \bar{a}_2,\bar{b}_2$ denote the standard generators of the fundamental group of a genus 2 Riemann surface satisfying $\prod_{i=1}^2[\bar{a}_i,\bar{b}_i]=1$.  Then we can assume that the inclusion $\bar{\Sigma}_2 \hookrightarrow T^4\#\CPb$ maps the generators of $\pi_1(\bar{\Sigma}_2)$ as follows:
\begin{equation}\label{eq: embedding of Sigma_2'}
\bar{a}_1  \mapsto  \alpha_1, \ \ 
\bar{b}_1  \mapsto  \alpha_2, \ \ 
\bar{a}_2  \mapsto  \alpha_3^2, \ \ 
\bar{b}_2  \mapsto  \alpha_4.
\end{equation}
Note that $\bar{a}_2$ is mapped to a square.

\section{Torus surgeries on $T^4\#\CPb$}
\label{sec: Luttinger surgery}

\begin{figure}[ht]
\begin{center}
\includegraphics[scale=.47]{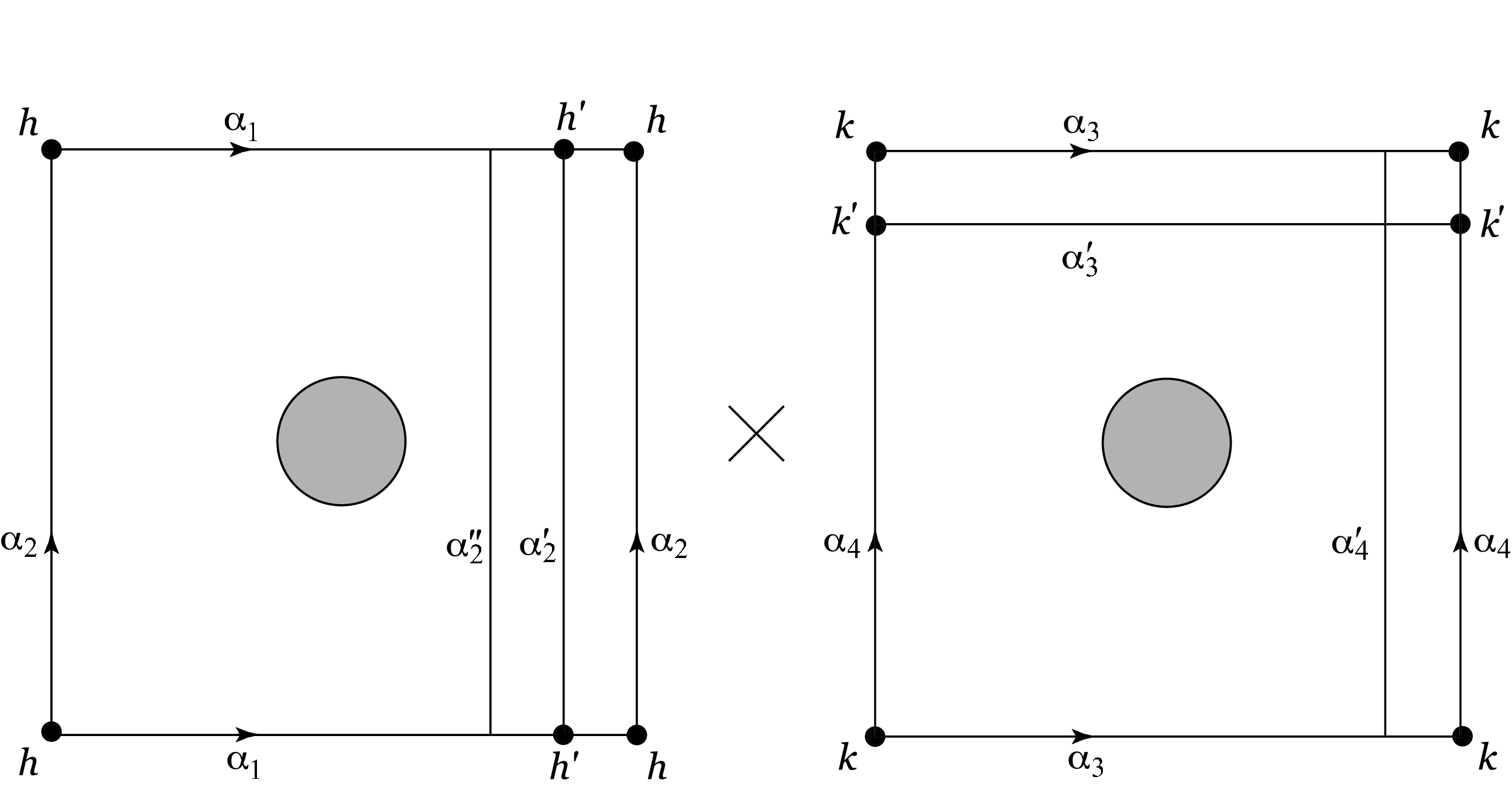}
\caption{Lagrangian tori $\alpha_2' \times \alpha_3'$ and $\alpha_2'' \times \alpha_4'$}
\label{fig:squares}
\end{center}
\end{figure}

Note that $\alpha_2\times\alpha_3$ and $\alpha_1\times\alpha_4$ are geometrically dual Lagrangian tori in $T^4\#\CPb$ and is disjoint from $\bar{\Sigma}_2$.  Let $Z'$ denote the symplectic 4-manifold that is the result of $(\alpha_2'\times\alpha_3',\alpha_3',-1)$ Luttinger surgery (cf.\ \cite{ADK, FPS}) on $T^4\#\CPb$.  The commutator relation $[\alpha_1,\alpha_4]=1$ in $\pi_1(T^4\#\CPb)\cong\pi_1(T^4)$ is replaced by the relation $\alpha_3=[\alpha_1^{-1},\alpha_4^{-1}]$ in $\pi_1(Z')$ (cf.\ \cite{BK:1-3, FPS}).  
We refer to Figure~\ref{fig:squares}, wherein we can see that the Lagrangian push-off of 
$\alpha_3'$ is $\alpha_3$ (also see Figure~1 in \cite{BK:1-3}).  
Note that $\bar{\Sigma}_2$ is still a symplectic submanifold of $Z'$.  

Now consider another pair of geometrically dual Lagrangian tori $\alpha_1\times\alpha_3$ and $\alpha_2\times\alpha_4$ in $T^4\#\CPb$.  Given any triple of integers $q\geq 0$, $r\geq 0$ and $m\geq 1$ satisfying $\gcd(m,r)=1$, let $Z''(1/q,m/r)$ denote the result of the following two torus surgeries on $T^4\#\CPb$:
\begin{equation}\label{eq: Luttinger surgeries in T^4}
(\alpha_2'\times\alpha_3',\alpha_3',-1/q), \ \ (\alpha_2''\times\alpha_4',\alpha_4',-m/r).  
\end{equation}
If $q=0$ or $r=0$, the corresponding surgery is trivial, i.e., the result of such surgery is diffeomorphic to the original $4$-manifold.  In particular, $Z''(1/0,1/0)=T^4\#\CPb$ and $Z''(1/1,1/0)=Z'$.  

Both surgeries in (\ref{eq: Luttinger surgeries in T^4}) are to be performed with respect to the Lagrangian framing and as such
the first is a Luttinger surgery for any $q\geq 1$ and the second is a Luttinger surgery when $m=1$ and $r\geq 1$.   In particular, $Z''(1/q,1/r)$ inherits a canonical symplectic structure from $T^4\#\CPb$.  Note that $\bar{\Sigma}_2$ is still a submanifold of $Z''(1/q,m/r)$ and a symplectic submanifold of $Z''(1/q,1/r)$.   

Finally, we note that $e(Z')=e(Z''(1/q,m/r))=1$ and $\sigma(Z')=\sigma(Z''(1/q,m/r))\linebreak[0] =-1$.  Neither $Z'$ nor $Z''(1/q,1/r)$ is a minimal symplectic $4$-manifold, but the pairs $(Z',\bar{\Sigma}_2)$ and $(Z''(1/q,1/r),\bar{\Sigma}_2)$ are both relatively minimal (cf.\ \cite{tjli}).  We should also point out that the torus surgeries in this section and in Section~\ref{sec:(2n-3)(S^2 x S^2)} are not the only ones we could have chosen.  Many other combinations of surgeries work just as well for the constructions that follow.

\section{Choice of base point}
\label{sec: base point}

\begin{figure}[ht]
\begin{center}
\includegraphics[scale=.49]{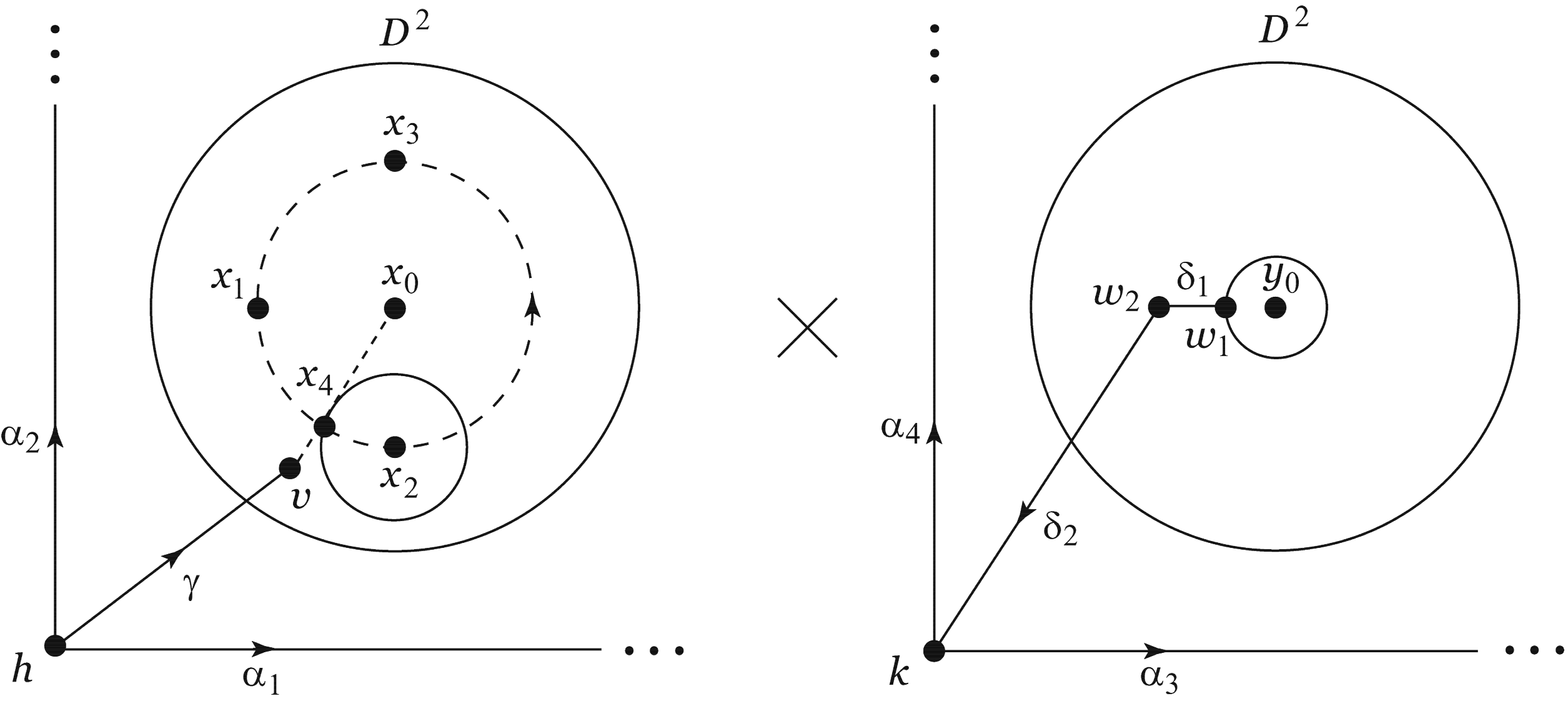}
\caption{Base point $h\times w_1$\/ and paths $\gamma$\/ and $\Delta=\delta_1\delta_2$}
\label{fig: disk}
\end{center}
\end{figure}

It will be convenient to 
view $T^2$ as the quotient of the unit square $[0,1]\times[0,1]=\R/\Z \times \R/\Z$ with the opposite edges identified.  Let $y_0=({1}/{2},{1}/{2})\in \alpha_3\times\alpha_4 = T^2$.  For $s\in [0,1]$, define
\begin{equation*}
\eta(s) = x_0+\|x_1-x_0\|e^{\pi i (2s+1)} \in \alpha_1\times \alpha_2 = T^2.
\end{equation*}
The point $\eta(s)$ is the counterclockwise rotation of the point $x_1$ about the center point $x_0=(1/2,1/2)\in D^2$ by angle $2\pi s$.  See Figure~\ref{fig: disk}.  The two large disks labeled by $D^2$ in Figure~\ref{fig: disk} are exactly the shaded disks in Figure~\ref{fig:squares}.  We identify the disk $D^2$ in the $\alpha_1\times\alpha_2$ torus with the embedded disk in $T^2$ that we introduced at the beginning of Section~\ref{sec:braided surface}.  It follows that the torus $T_{\beta}$ in Section~\ref{sec:braided surface} is smoothly isotopic to the symplectic torus in $T^4$ parameterized by 
\begin{equation}\label{eq:T_beta}
(s,t) \,\mapsto\, \eta(s) \times (2s, t) \,\in\, D^2 \times T^2 \,\subset\, T^2\times T^2 ,
\end{equation}
where $(s,t)\in [0,1]\times[0,1]$.  
By looking at Figure~\ref{fig:squares}, we can immediately deduce the following.  

\begin{lem}\label{lem: subsquares}
Let $R$ be the radius of the shaded disks $D^2$ in Figure\/ $\ref{fig:squares}$.  Let\/
$[0,\frac{1}{2}+R]^2$ denote the cartesian product that is a lower left sub-square of the unit square $[0,1]^2$ containing $D^2$.  Then the union 
\begin{equation}\label{eq: subsquares}
\big( T^2 \times [0, {\textstyle \frac{1}{2}} +R]^2 \big) \,\cup\, 
\big( [0, {\textstyle \frac{1}{2}} +R]^2 \times T^2 \big)
\end{equation}
is disjoint from the two surgery tori $\alpha_2' \times\alpha_3'$ and $\alpha_2''\times\alpha_4'$.  Thus, away from the blow-up point, $(\ref{eq: subsquares})$ can be viewed as a subset of $Z''(1/q,m/r)$.  
\hfill $\qed$
\end{lem}

Next we note that the intersection between $(\alpha_1\times\alpha_2)\times y_0=T^2\times y_0$ and $T_{\beta}$ consists of two points $x_2 \times y_0$ and $x_3 \times y_0$, which correspond to $s=1/4$ and $s=3/4$.  Our genus 2 surface $\bar{\Sigma}_2$ is obtained from the union of $T^2\times y_0$ and $T_{\beta}$ by resolving the intersection $x_2 \times y_0$ and blowing up at $x_3\times y_0$.  We will assume that the resolution and the blow-up take place inside open 4-dimensional balls of fixed radius $\rho$\/ centered at $x_2\times y_0$ and $x_3\times y_0$, respectively.  This means that the points on $T^2\times y_0$ and $T_{\beta}$ that lie outside these open 4-balls are also points on $\bar{\Sigma}_2$.  

Let $\nu\bar{\Sigma}_2$ be the $\epsilon$-tubular neighborhood of $\bar{\Sigma}_2$, i.e., an open subset of $Z''(1/q,m/r)$ consisting of points that are within distance $\epsilon>0$ from $\bar{\Sigma}_2$.  Here, we require $\epsilon$\/ to be much smaller than $\rho$, which in turn is much smaller than $R$.  
We need $\epsilon \ll \rho$\/ so that the parts of $\nu\bar{\Sigma}_2$ inside the open 4-balls of radius $\rho$\/ centered at $x_2\times y_0$ and $x_3\times y_0$ will look like the standard local models. 
We need $\rho \ll R$\/ so that the aforementioned 4-balls are both thoroughly contained in $D^2 \times D^2$.  

Let $D^2_{\epsilon}(y_0)$ be the open disk of radius $\epsilon$\/ in the $\alpha_3\times\alpha_4$ torus centered at $y_0$ in Figure~\ref{fig: disk}.  (The scales in Figure~\ref{fig: disk} are exaggerated quite a bit for the ease of labeling.  The disk $D^2_{\epsilon}(y_0)$ is quite small.)  Let $w_1=(\frac{1}{2}-\epsilon,\frac{1}{2})\in\partial D^2_{\epsilon}(y_0)$.  Note that $(s,t) \times D^2_{\epsilon}(y_0)$ is a normal disk to $\bar{\Sigma}_2$ at the point 
\begin{equation*}
(s,t)\times y_0 \,\in\, \big[T^2 \setminus \big(D^2_{\rho}(x_2) \cup D^2_{\rho}(x_3)\big)\big]\times y_0 \,\subset\, \bar{\Sigma}_2
\end{equation*} 
for any $(s,t)\in \alpha_1\times\alpha_2$ that lies outside the two open disks of radius $\rho$\/ centered at $x_2$ and $x_3$ wherein the resolution and the blow-up take place.  The disk $D^2_{\rho}(x_2)$ of radius $\rho$\/ centered at $x_2$ in the $\alpha_1\times\alpha_2$ torus is drawn in Figure~\ref{fig: disk}.  We choose small enough $\rho$\/ so that $D^2_{\rho}(x_2)$ and $D^2_{\rho}(x_3)$ are both contained in $D^2 \setminus\{x_0\}$.  

Let $x_4$ be the point shown in Figure~\ref{fig: disk} that lies in the $\alpha_1\times\alpha_2$ torus at distance $\|x_1-x_0\|$ from $x_0$ and at distance $\rho$\/ from $x_2$.  We have $x_4=\eta(s_4)$ for some $0 < s_4 < 1/4$.  From (\ref{eq:T_beta}) and our choice of $x_4$, we know that $x_4\times (2s_4,t)$ is a point in $\bar{\Sigma}_2$ for all $t\in [0,1]$.  Let $w_2=(2s_4,1/2)$, and let $v$\/ be a point in the $\alpha_1\times\alpha_2$ torus in Figure~\ref{fig: disk} that lies on the line through $x_0$ and $x_4$, and is at distance $\epsilon$\/ from $x_4$.  Note that $v\times w_2$ lies on the boundary circle of the normal disk to $\bar{\Sigma}_2$ of radius $\epsilon$\/ centered at $x_4\times w_2\in \bar{\Sigma}_2$.  

Let $\delta_1$ be the horizontal line path in the $\alpha_3\times\alpha_4$ torus from $w_1$ to $w_2$.  Let $\delta_2$ be the straight line path from $w_2$ to the corner point $k$.  Both paths lie away from the two surgery tori $\alpha_2' \times\alpha_3'$ and $\alpha_2''\times\alpha_4'$.  Let $\Delta=\delta_1\delta_2$ denote the piecewise linear path that is the composition of these two line paths: $\Delta(t)=\delta_1(2t)$ for $t\in [0,1/2]$ and $\Delta(t)=\delta_2(2t-1)$ for $t\in [1/2,1]$.  (We always write path composition from left to right.)  

In applications that follow, the base point of our fundamental group needs to lie on the boundary in order to apply Seifert-Van Kampen theorem.  
We choose the base point of $\pi_1(Z''(1/q,m/r)\setminus\nu\bar{\Sigma}_2)$ to be the point $h\times w_1$ in Figure~\ref{fig: disk}.  Note that $h\times w_1$ lies on the boundary circle of a normal disk $h \times D^2_{\epsilon}(y_0)$ to $\bar{\Sigma}_2$ at the point $h\times y_0\in\bar{\Sigma}_2$.  Hence our base point $h\times w_1$\/ is indeed on the boundary of the $\epsilon$-tubular neighborhood of $\bar{\Sigma}_2$.  Later we will choose our parallel copy $\bar{\Sigma}_2^{\scparallel}$ so that it contains $h\times w_1$.

\section{Inclusion induced homomorphism}
\label{sec: inclusion}

The generator $\alpha_1$ in $\pi_1(Z''(1/q,m/r)\setminus\nu\bar{\Sigma}_2, h\times w_1)$ is represented by a based loop that is the composition of three paths: $h\times\Delta$ from $h\times w_1$ to $h\times k$, the parameterized loop $t\mapsto (t,0,0,0)$ for $t\in[0,1]$, and the reverse path $h\times\Delta^{-1}$ from $h\times k$\/ to $h\times w_1$, in that order.  Similarly, the generators $\alpha_i$, $i=2,3,4$, can be represented by loops based at $h\times w_1$ that are conjugates of loops lying on the boundaries of the squares by the fixed path $h\times\Delta$.  The images of these generators in $\pi_1(Z''(1/q,m/r), h\times w_1)$ under the homomorphism induced by the inclusion $Z''(1/q,m/r)\setminus\nu\bar{\Sigma}_2 \hookrightarrow Z''(1/q,m/r)$ will continue to be denoted by $\alpha_1,\dots,\alpha_4$.  

Let $\partial(\nu\bar{\Sigma}_2)$ denote the boundary of the $\epsilon$-tubular neighborhood of $\bar{\Sigma}_2$ inside $Z''(1/q,m/r)$. 
We fix a normal projection map ${\rm pr}:\partial(\nu\bar{\Sigma}_2)\rightarrow \bar{\Sigma}_2$ such that we have ${\rm pr}(h\times w_1)=h\times y_0$.  

\begin{lem}\label{lem: inclusion}
The inclusion\/ $\bar{\Sigma}_2\hookrightarrow Z''(1/q,m/r)$ maps the standard generators $\bar{a}_i,\bar{b}_i$ $(i=1,2)$ of\/ $\pi_1(\bar{\Sigma}_2, h\times y_0)$ as in\/ $(\ref{eq: embedding of Sigma_2'})$.  
\end{lem}

\begin{proof}
Consider the following simple parameterized loops in (\ref{eq: subsquares}) that are based at $x_2\times y_0$ with $t\in[0,1]$:
\begin{equation}\label{eq: bar-generator}
\begin{array}{l}
\bar{a}_1(t) = x_2\times y_0 + (t,0,0,0),  \\
\bar{b}_1(t) = x_2\times y_0 + (0,t,0,0),  \\ 
\bar{a}_2(t) = (x_0+\|x_1-x_0\|e^{\pi i (2t+\frac{3}{2})}) \times (2t+\frac{1}{2},\frac{1}{2}), \\
\bar{b}_2(t) = x_2\times y_0 + (0,0,0,t). 
\end{array}
\end{equation}
Here, the coordinates are in $\R/\Z$.  
Note that $\bar{a}_1(t)$ and $\bar{b}_1(t)$ lie in $T^2 \times y_0$, and they represent the standard generators of $\pi_1(T^2\times y_0,x_2\times y_0)$.  Similarly, $\bar{a}_2(t)$ and $\bar{b}_2(t)$ lie in $T_{\beta}$ and represent the standard generators of $\pi_1(T_{\beta},x_2\times y_0)$.
We will describe how the above four loops uniquely determine four standard generators of $\pi_1(\bar{\Sigma}_2, h\times y_0)$.  

By Lemma~\ref{lem: subsquares}, 
we observe that all four loops are disjoint from the two surgery tori $\alpha_2' \times\alpha_3'$ and $\alpha_2''\times\alpha_4'$, and thus they descend to simple curves in $Z''(1/q,m/r)$ away from the blow-up point $x_3\times y_0$.  Note that $\bar{a}_1(t)$ and $\bar{b}_2(t)$ avoid $x_3\times y_0$, but $\bar{b}_1(t)$ and $\bar{a}_2(t)$ intersect each other transversely at  $x_3\times y_0$.  There is a canonical isomorphism $\pi_1(T^4\#\CPb)\cong\pi_1(T^4)$, so the above four based loops lift uniquely (up to homotopy) to four loops in $Z''(1/q,m/r)$ that intersect each other only at the point $x_2\times y_0$.  
To obtain $\bar{\Sigma}_2$, we need to resolve the intersection point $x_2\times y_0$ between $T^2\times y_0$ and $T_{\beta}$.  The resolution procedure takes place inside a 4-ball of a small radius $\rho$\/ centered at $x_2\times y_0$ and is completely local in nature.  In particular, the resolution does not change the fundamental group of the ambient 4-manifold.  

Let $\Gamma$ be the straight line path from $h\times y_0$ to $x_2 \times y_0$ in $Z''(1/q,m/r)$.  After conjugating by $\Gamma$, $\bar{a}_i(t),\bar{b}_i(t)$ ($i=1,2$) give rise to loops based at $h\times y_0$, which we continue to denote by $\bar{a}_i(t),\bar{b}_i(t)$ for simplicity.  Clearly, 
there are based homotopies inside $Z''(1/q,m/r)$ from $\bar{a}_i(t),\bar{b}_i(t)$ ($i=1,2$) to four standard generators of $\pi_1(\bar{\Sigma}_2, h\times y_0)$.  In fact, we can arrange that such homotopies are supported completely inside 
(\ref{eq: subsquares}).  
Let $\bar{a}_i,\bar{b}_i$ ($i=1,2$) denote the generators of $\pi_1(\bar{\Sigma}_2, h\times y_0)$ corresponding to loops $\bar{a}_i(t),\bar{b}_i(t)$, respectively.  

Let $\Phi$ denote the straight line path in the closure of $h\times D^2_{\epsilon}(y_0)$ from $h\times w_1$ to ${\rm pr}(h\times w_1)=h\times y_0$.  After conjugating by $\Phi$, $\bar{a}_i(t),\bar{b}_i(t)$ ($i=1,2$) can be viewed as loops based at $h\times w_1$.  (I.e., we conjugate the parameterized loops in (\ref{eq: bar-generator}) by the composition of paths $\Phi\Gamma$.)  

We first consider $\bar{a}_2(t)$, which is the most difficult case.  
Since the blow-up process preserves $\pi_1$ canonically, we can smoothly perturb away the part of $\bar{a}_2(t)$ lying in $\CPb\setminus D^4_{\rho}$ (which is simply-connected) such that the image of the new $\bar{a}_2(t)$ lies inside (\ref{eq: subsquares}) minus a 4-ball of radius $\rho$\/ centered at the blow-up point $x_3\times y_0$.  Near $x_3\times y_0$, we require the second and the fourth coordinates of the new $\bar{a}_2(t)$ to be smaller than those of the old $\bar{a}_2(t)$.  More descriptively, we are nudging the loop $\bar{a}_2(t)$ slightly downward near $x_3\times y_0$ to avoid the blow-up area so that we may write down our formulas using the coordinates inherited from $T^4$. 
If $\bar{a}_2(t)=(\xi_1(t),\xi_2(t),\xi_3(t),\xi_4(t))\in [0,\frac{1}{2}+R]^2 \times [0,1]^2$, for $t\in [0,1]$, denotes this new parameterized loop based at $h\times w_1$, then we define 
\begin{equation*}
G(s,t)=\big((1-s)\xi_1(t),(1-s)\xi_2(t),\xi_3(t),(1-s)\xi_4(t)\big).
\end{equation*}
A based homotopy $F_3:[0,1]\times [0,1]\rightarrow Z''(1/q,m/r)$ from $\bar{a}_2(t)$ to $\alpha_3^2$ is now given by 
\begin{equation*}
F_3(s,t) = \left\{ \begin{array}{ccl} 
\big(0,0,\frac{1}{2}-\epsilon,\frac{1-3t}{2}\big) & \text{ if } & 0\leq t \leq \frac{s}{3}, \\[2pt]
G\big(s ,\,\frac{3t-s}{3-2s}\big) & \text{ if } & \frac{s}{3} \leq t \leq 1-\frac{s}{3}, \\[2pt] 
\big(0,0,\frac{1}{2}-\epsilon,\frac{3t-2}{2}\big) & \text{ if } & 1-\frac{s}{3}\leq t \leq 1.  
\end{array}\right.
\end{equation*}

More descriptively, $F_3$ contracts the first and the second coordinates of $\bar{a}_2(t)$ to the lower left corner point $h$\/ of the left $T^2$ in Figure \ref{fig: disk}, and projects the third and the fourth coordinates of $\bar{a}_2(t)$ vertically down onto the third axis of $T^4$.  
We can check that the first and the second coordinates of $F_3$ are contained in the lower left sub-square $[0,\frac{1}{2}+R]$.  
Hence the image of $F_3$ is contained in (\ref{eq: subsquares}), and is disjoint from the two surgery tori by Lemma~\ref{lem: subsquares}.  

Similarly, $\bar{a}_1(t)$, $\bar{b}_1(t)$ and $\bar{b}_2(t)$ are based homotopic to $\alpha_1$, $\alpha_2$ and $\alpha_4$, respectively.  Once again, these based homotopies can be supported completely in (\ref{eq: subsquares}).  
\end{proof}

\section{Parallel loops on the boundary}
\label{sec: parallel}

Let $\bar{\Sigma}_2^{\scparallel}$ denote a parallel copy of $\bar{\Sigma}_2$ in the boundary of $\epsilon$-tubular neighborhood $\partial (\nu\bar{\Sigma}_2) \subset Z''(1/q,m/r)$ through the base point $h\times w_1$.  Let $\bar{a}_i^{\scparallel}$, $\bar{b}_i^{\scparallel}$, $i=1,2$, be the generators of $\pi_1(\bar{\Sigma}_2^{\scparallel})$ that are push-offs of the generators $\bar{a}_i$, $\bar{b}_i$ of $\pi_1(\bar{\Sigma}_2)$.  

\begin{lem}\label{lem: parallel}
There exists an embedding\/ $\bar{\Sigma}_2^{\scparallel}\hookrightarrow Z''(1/q,m/r)\setminus\nu\bar{\Sigma}_2$ such that the induced 
homomorphism on\/ $\pi_1$ maps the generators\/ $\bar{a}_1^{\scparallel}$, $\bar{b}_1^{\scparallel}$ and\/ $\bar{b}_2^{\scparallel}$ of\/ $\pi_1(\bar{\Sigma}_2^{\scparallel})$ to\/ $\alpha_1$, $\alpha_2$ and\/ $\alpha_4$ in\/ $\pi_1(Z''(1/q,m/r)\setminus\nu\bar{\Sigma}_2)$, respectively.  
\end{lem}

\begin{proof}
We claim that there exist parameterized loops $\bar{a}_1^{\scparallel}(t)$, $\bar{b}_1^{\scparallel}(t)$ and $\bar{b}_2^{\scparallel}(t)$ in $\partial (\nu\bar{\Sigma}_2)$ that are based at $h\times w_1$ and represent $\alpha_1$, $\alpha_2$ and $\alpha_4$ respectively in $\pi_1(Z''(1/q,m/r)\setminus\nu\bar{\Sigma}_2,h\times w_1)$.  Our claim implies the lemma as follows.  

First we fix a diffeomorphism $\phi : \partial (\nu\bar{\Sigma}_2) \rightarrow S^1\times \bar{\Sigma}_2$.  Let ${\rm pr}_1 : S^1\times \bar{\Sigma}_2 \rightarrow S^1$ and ${\rm pr}_2 : S^1\times \bar{\Sigma}_2 \rightarrow \bar{\Sigma}_2$ be the projections onto the factors.  Choose an orientation of $S^1$ such that $\phi$\/ is orientation preserving, and let $\zeta$\/ denote the corresponding generator of $H^1(S^1;\Z)\cong \Z$.  

Any parallel copy $\bar{\Sigma}_2^{\scparallel}$ is the image of some push-off embedding $f : \bar{\Sigma}_2 \rightarrow \partial (\nu\bar{\Sigma}_2)$ 
such that the composition ${\rm pr}_2 \circ \phi \circ f : \bar{\Sigma}_2 \rightarrow \bar{\Sigma}_2$ is isotopic to the identity map.  Let 
$\bar{a}_i^{\scparallel}=f_{\ast}(\bar{a}_i)$ and $\bar{b}_i^{\scparallel}=f_{\ast}(\bar{b}_i)$ ($i=1,2$).  The pullback of $\zeta$\/ under the composition 
\begin{equation*}
\bar{\Sigma}_2 \stackrel{f}{\longrightarrow} \partial (\nu\bar{\Sigma}_2) 
\stackrel{\phi}{\longrightarrow} S^1\times \bar{\Sigma}_2 
\stackrel{{\rm pr}_1}{\longrightarrow} S^1
\end{equation*}
gives an element $\tau_f \in H^1(\bar{\Sigma}_2; \Z)$.  
Conversely, any element of $H^1(\bar{\Sigma}_2; \Z)$ can be expressed as $\theta^{\ast}(\zeta)$ for some smooth map $\theta : \bar{\Sigma}_2 \rightarrow S^1$.  The image of the graph of $\theta$\/ under $\phi^{-1}$, 
\begin{equation*}
\{ \phi^{-1}(\theta(x),x) \in \partial (\nu\bar{\Sigma}_2) \mid x \in \bar{\Sigma}_2 \},
\end{equation*}
is a parallel copy $\bar{\Sigma}_2^{\scparallel}$.  
It is well known that this gives a 
one-to-one correspondence between the isotopy classes of $\bar{\Sigma}_2^{\scparallel}$ and the elements of $H^1(\bar{\Sigma}_2; \Z)$.  

Next, $\phi$\/ gives rise to a factorization 
\begin{equation}\label{eq: pi_1 factorization}
\pi_1(\partial (\nu\bar{\Sigma}_2)) \,\cong\, 
\pi_1(S^1) \times \pi_1(\bar{\Sigma}_2) \,\cong\, 
\Z \times \pi_1(\bar{\Sigma}_2).
\end{equation}
The $\Z$ factor in (\ref{eq: pi_1 factorization}) is generated by a meridian $\mu(\bar{\Sigma}_2)$ of $\bar{\Sigma}_2$, and every element in the $\Z$ factor commutes with any element in $\pi_1(\partial (\nu\bar{\Sigma}_2))$.  By Seifert-Van Kampen theorem, we have
\begin{equation}\label{eq: fill in}
\pi_1(Z''(1/q,m/r)) \cong \frac{\pi_1(Z''(1/q,m/r)\setminus\nu\bar{\Sigma}_2) 
\ast \pi_1(\nu\bar{\Sigma}_2)}
{\langle \mu(\bar{\Sigma}_2)=1,\, \bar{a}_i^{\scparallel}=\bar{a}_i,\, \bar{b}_i^{\scparallel}=\bar{b}_i \ (i=1,2)\rangle}.
\end{equation}
Note that when we glue $Z''(1/q,m/r)\setminus\nu\bar{\Sigma}_2$ and $\nu\bar{\Sigma}_2$ together, the resulting closed $4$-manifold is always diffeomorphic to $Z''(1/q,m/r)$ for every possible choice of $\bar{\Sigma}_2^{\scparallel}$ and every boundary identification 
(cf.\ Remark 8.1.3.(a) in \cite{GS}).

By Lemma~\ref{lem: inclusion}, the loops $\bar{a}_1^{\scparallel}(t)$, $\bar{b}_1^{\scparallel}(t)$ and $\bar{b}_2^{\scparallel}(t)$ in $\partial (\nu\bar{\Sigma}_2)=\partial(Z''(1/q,m/r)\setminus\nu\bar{\Sigma}_2)$ representing $\alpha_1$, $\alpha_2$ and $\alpha_4$ in $\pi_1(Z''(1/q,m/r)\setminus\nu\bar{\Sigma}_2)$ determine the generators $\bar{a}_1$, $\bar{b}_1$ and $\bar{b}_2$ of $\pi_1(\bar{\Sigma}_2)$, respectively, via (\ref{eq: fill in}).  Therefore, we have 
\begin{equation*}
({\rm pr}_2 \circ \phi)_{\ast}(\bar{a}_1^{\scparallel}(t)) = \bar{a}_1, \ \ 
({\rm pr}_2 \circ \phi)_{\ast}(\bar{b}_1^{\scparallel}(t)) = \bar{b}_1, \ \
({\rm pr}_2 \circ \phi)_{\ast}(\bar{b}_2^{\scparallel}(t)) = \bar{b}_2,
\end{equation*} 
where $({\rm pr}_2 \circ \phi)_{\ast}: \pi_1(\partial (\nu\bar{\Sigma}_2)) \rightarrow \pi_1(\bar{\Sigma}_2)$ is the induced homomorphism.

Finally, we choose $\bar{\Sigma}_2^{\scparallel}=f(\bar{\Sigma}_2)$ such that the corresponding $\tau_f \in H^1(\bar{\Sigma}_2; \Z)\cong {\rm Hom}(H_1(\bar{\Sigma}_2;\Z),\Z)$ satisfies 
\begin{equation*}
\tau_f(\bar{a}_1)=({\rm pr}_1 \circ \phi)_{\ast} (\bar{a}_1^{\scparallel}(t)), \  
\tau_f(\bar{b}_1)=({\rm pr}_1 \circ \phi)_{\ast} (\bar{b}_1^{\scparallel}(t)), \  
\tau_f(\bar{b}_2)=({\rm pr}_1 \circ \phi)_{\ast} (\bar{b}_2^{\scparallel}(t)),
\end{equation*}
where $({\rm pr}_1 \circ \phi)_{\ast}: \pi_1(\partial (\nu\bar{\Sigma}_2)) \rightarrow \pi_1(S^1) \cong \Z$.  From (\ref{eq: pi_1 factorization}), it now follows that $\bar{a}_1^{\scparallel}(t)$, $\bar{b}_1^{\scparallel}(t)$ and $\bar{b}_2^{\scparallel}(t)$ represent $\bar{a}_1^{\scparallel}$, $\bar{b}_1^{\scparallel}$ and\/ $\bar{b}_2^{\scparallel}$ in $\pi_1(\partial (\nu\bar{\Sigma}_2))$, respectively.  
This proves the lemma.  

We should note that it is not necessary to check whether the three loops $\bar{a}_1^{\scparallel}(t)$, $\bar{b}_1^{\scparallel}(t)$ and $\bar{b}_2^{\scparallel}(t)$ in $\partial (\nu\bar{\Sigma}_2)$ are completely contained in some particular $\bar{\Sigma}_2^{\scparallel}$ since it suffices to know the fundamental group elements they represent in (\ref{eq: pi_1 factorization}).  In other words, the existence of $\bar{a}_1^{\scparallel}(t)$, $\bar{b}_1^{\scparallel}(t)$ and $\bar{b}_2^{\scparallel}(t)$ in our claim implies the existence of three loops in some $\bar{\Sigma}_2^{\scparallel}$ representing the same fundamental group elements, and for the lemma we only need to pin down the fundamental group elements.  

It only remains to exhibit the three loops in our initial claim.  
For $t\in[0,1]$, let $\bar{a}_1^{\scparallel}(t)=(t,0,\frac{1}{2}-\epsilon,\frac{1}{2})=h\times w_1 + (t,0,0,0)$ be a simple loop based at $h\times w_1$.  For each $t\in [0,1]$, the corresponding point $\bar{a}_1^{\scparallel}(t)$ is in $\partial(\nu\bar{\Sigma}_2)$ since it lies on the boundary circle of a normal disk $(t,0) \times D^2_{\epsilon}(y_0)$ to $\bar{\Sigma}_2$ at the point $(t,0) \times y_0 \in T^2 \times y_0$.  Each point $(t,0)\times y_0$ lies in $\bar{\Sigma}_2$ since the point $(t,0)$ in the $\alpha_1\times\alpha_2$ torus lies at distance more than $\rho$\/ from the intersection coordinates $x_2$ and $x_3$.  Similarly, let $\bar{b}_1^{\scparallel}(t)=(0,t,\frac{1}{2}-\epsilon,\frac{1}{2})=h\times w_1 + (0,t,0,0)$ be another simple loop based at $h\times w_1$.  For each $t\in [0,1]$, the image point $\bar{b}_1^{\scparallel}(t)\in\partial(\nu\bar{\Sigma}_2)$ lies on the boundary circle of a normal disk $(0,t) \times D^2_{\epsilon}(y_0)$ to $\bar{\Sigma}_2$ at !
 the point $(0,t)\times y_0\in \bar{\Sigma}_2$.  

Next, let $\gamma$\/ denote the straight line path in the $\alpha_1\times\alpha_2$ torus from the corner point $h$\/ to the point $v$\/ in Figure~\ref{fig: disk}.  
Let $\delta_1$ be the horizontal line path in the $\alpha_3\times\alpha_4$ torus from $w_1$ to $w_2=(2s_4,1/2)$ as in Section~\ref{sec: base point}.  
We define the loop $\bar{b}_2^{\scparallel}(t)$, based at $h\times w_1$, to be the composition of five paths:  $\gamma\times w_1$, $v\times\delta_1$, the parameterized loop $t\mapsto v\times w_2 + (0,0,0,t)$ for $t\in[0,1]$, $v\times\delta_1^{-1}$, and $\gamma^{-1}\times w_1$, in that order.  Each point on the path $\gamma\times w_1$ (and its reverse $\gamma^{-1}\times w_1$) is in $\partial(\nu\bar{\Sigma}_2)$ since it lies on the boundary circle of a normal disk to $\bar{\Sigma}_2$ at some point $\gamma(t)\times y_0 \in\bar{\Sigma}_2$.  

It is also not hard to check that each point in the middle loop, $v\times w_2 + (0,0,0,t)$, is in $\partial(\nu\bar{\Sigma}_2)$.  In Section~\ref{sec: base point}, we have already observed that $x_4\times (2s_4,t)$ lies in the $T_{\beta}$ part of $\bar{\Sigma}_2$ for all $t\in [0,1]$.  Since $v\times w_2$ is on the boundary circle of a normal disk to $\bar{\Sigma}_2$ at $x_4\times (2s_4,1/2)$, we see that each point $v\times w_2 + (0,0,0,t)$ lies on the boundary circle of a normal disk to $\bar{\Sigma}_2$ at $x_4\times (2s_4,\frac{1}{2}+t)$. 

The remaining path $v\times \delta_1$ (and its reverse $v\times \delta_1^{-1}$) does not lie in $\partial(\nu\bar{\Sigma}_2)$.  However, this path lies in the complement $Z''(1/q,m/r)\setminus\nu\bar{\Sigma}_2$ and the endpoints of this path, $v\times w_1$ and $v\times w_2$, are both in $\partial(\nu\bar{\Sigma}_2)$.  Moreover, the straight line path $v\times \delta_1$ is the \emph{shortest}\/ path in $Z''(1/q,m/r)\setminus\nu\bar{\Sigma}_2$ between $v\times w_1$ and $v\times w_2$.  It is not hard to see that there will be a path $\delta_1^{\scparallel}$ from $v\times w_1$ to $v\times w_2$ that lies completely in the boundary $\partial(\nu\bar{\Sigma}_2)$ such that $\delta_1^{\scparallel}$ is homotopic to $v\times \delta_1$ inside the complement $Z''(1/q,m/r)\setminus\nu\bar{\Sigma}_2$ fixing the endpoints.  Such homotopy can be supported inside (\ref{eq: subsquares}), away from the two surgery tori $\alpha_2' \times\alpha_3'$ and $\alpha_2''\times\alpha_4'$.  

Recall that when we resolve the intersection point $x_2\times y_0$, we replace two transversely intersecting disks near $x_2\times y_0$ with an unknotted cylinder $C\subset \bar{\Sigma}_2$ that connects the boundary circles of these two disks.  The path $\delta_1^{\scparallel}$ lies inside the path-connected set  
\begin{equation}\label{eq: nu C}
\partial(\nu C) \cap \partial(\nu\bar{\Sigma}_2) \,\cong\,
S^1 \times C \,\cong\, S^1 \times \big(S^1 \times [0,1]\big)
\,\cong\, T^2\times [0,1],
\end{equation}  
and it does not wind around any homotopically nontrivial circles in (\ref{eq: nu C}).  
This is because $v\times\delta_1$ is length-minimizing and the arc-length is a continuous function over paths that are homotopic to $v\times\delta_1$.  We choose $\delta_1^{\scparallel}$ such that it is length-minimizing among all the paths in $\partial(\nu C)\cap \partial(\nu\bar{\Sigma}_2)$ that start at $v\times w_1$ in the boundary component corresponding to $T^2 \times 0$ and end at $v\times w_2$ in the boundary component corresponding to $T^2 \times 1$.  When we project $\delta_1^{\scparallel}$ to different factors in (\ref{eq: nu C}), it will surject onto the $[0,1]$ factor, but it will have very small footprint on the $T^2$ factor.  Paths with the same endpoints as $\delta_1^{\scparallel}$ that wind around nontrivially in the $T^2$ factor cannot be length-minimizing since they will have more arc-length contribution from the $T^2$ factor.  See Figure~\ref{fig: cylinder} for a schematic diagram.  

\begin{figure}[ht]
\begin{center}
\includegraphics[scale=.65]{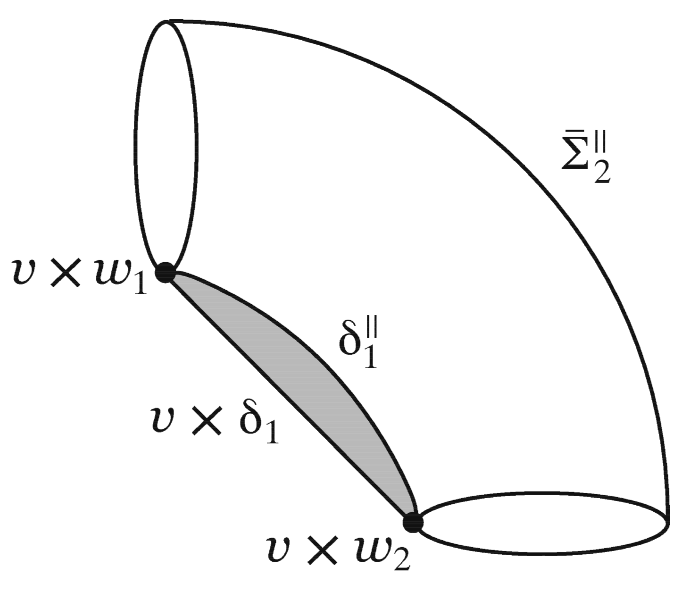}
\caption{Homotopy between $\delta_1^{\scparallel}$ and $v\times \delta_1$}
\label{fig: cylinder}
\end{center}
\end{figure}

In conclusion, $\bar{b}_2^{\scparallel}(t)$ does not lie completely in $\partial(\nu\bar{\Sigma}_2)$, but it is based homotopic to a loop in $\partial(\nu\bar{\Sigma}_2)$ (which is obtained by replacing $v\times\delta_1$ and its reverse with $\delta_1^{\scparallel}$ and its reverse, respectively) and such homotopy is supported in the complement $Z''(1/q,m/r)\setminus\nu\bar{\Sigma}_2$.  Since we are interested in the fundamental group of $Z''(1/q,m/r)\setminus\nu\bar{\Sigma}_2$, it will suffice to work with $\bar{b}_2^{\scparallel}(t)$.  

A based homotopy $F_1: [0,1]\times [0,1]\rightarrow Z''(1/q,m/r)\setminus\nu\bar{\Sigma}_2$ from $\bar{a}_1^{\scparallel}(t)$ to $\alpha_1$ and a based homotopy $F_2: [0,1]\times [0,1]\rightarrow Z''(1/q,m/r)\setminus\nu\bar{\Sigma}_2$ from $\bar{b}_1^{\scparallel}(t)$ to $\alpha_2$ are given respectively by 
\begin{equation*}
F_1(s,t) = \left\{ \begin{array}{ccl} 
h \times \Delta(3t) & \text{ if } &  0\leq t \leq \frac{s}{3}, \\[2pt]
\big( \frac{3t-s}{3-2s},0 \big) \times \Delta(s) & \text{ if } & \frac{s}{3} \leq t \leq 1-\frac{s}{3}, \\[2pt]  
h \times \Delta(3-3t) & \text{ if } & 1-\frac{s}{3}\leq t \leq 1,  
\end{array}\right. 
\end{equation*}
and
\begin{equation*}
F_2(s,t) = \left\{ \begin{array}{ccl} 
h \times \Delta(3t) & \text{ if } & 0\leq t \leq \frac{s}{3}, \\[2pt]
\big( 0, \frac{3t-s}{3-2s} \big) \times \Delta(s) & \text{ if } & \frac{s}{3} \leq t \leq 1-\frac{s}{3}, \\[2pt]  
h \times \Delta(3-3t) & \text{ if } & 1-\frac{s}{3}\leq t \leq 1.  
\end{array}\right. 
\end{equation*}
Similarly, a based homotopy $F_4: [0,1]\times [0,1]\rightarrow Z''(1/q,m/r)\setminus\nu\bar{\Sigma}_2$ from $\bar{b}_2^{\scparallel}(t)$ to $\alpha_4$ is given by 
\begin{equation*}
F_4(s,t) = \left\{ \begin{array}{ccl} 
\gamma(5t) \times \Delta(\frac{s(1-s)}{10}) & \text{ if } & 0\leq t \leq \frac{1-s}{5}, \\[2pt]
\gamma(1-s) \times \Delta \big(\frac{5}{2}(t-\frac{1}{5})+s(\frac{t}{2}+\frac{1}{2})\big) & \text{ if } & \frac{1-s}{5}\leq t \leq \frac{6-s}{15}, \\[2pt]
\gamma(1-s)\times \big[\Delta(\frac{1}{2}+\frac{8s}{15}-\frac{s^2}{30})+\big( 0, \frac{15t+s-6}{2s+3} \big)\big] & \text{ if } & \frac{6-s}{15} \leq t \leq \frac{9+s}{15}, \\[2pt]
\gamma(1-s)\times \Delta\big(\frac{5}{2}(\frac{4}{5}-t)+s(1-\frac{t}{2})\big) & \text{ if } & \frac{9+s}{15} \leq t \leq \frac{4+s}{5}, \\[2pt]
\gamma (5-5t) \times \Delta(\frac{s(1-s)}{10}) & \text{ if } & \frac{4+s}{5}\leq t \leq 1.  
\end{array}\right. 
\end{equation*}

We can check that the image of each $F_i$ ($i=1,2,4$) lies in (\ref{eq: subsquares}) and hence away from the two surgery tori $\alpha_2' \times\alpha_3'$ and $\alpha_2''\times\alpha_4'$.  More descriptively, each $F_i$ adds a non-simple part to each of the loops $\bar{a}_1^{\scparallel}(t)$, $\bar{b}_1^{\scparallel}(t)$ and $\bar{b}_2^{\scparallel}(t)$ in the third and fourth coordinates.  For fixed $s$, the parameterized loop $t \mapsto F_i(s,t)$ traverses down the path $\Delta$ in the third and fourth coordinates, makes a simple loop in the $i$-th direction, and then comes back up along $\Delta$.  The new non-simple part along $\Delta$ gets incrementally bigger as $s$\/ increases, eventually going all the way down to the lower left corner point $k=\Delta(1)$ of the right square in Figure~\ref{fig: disk}.  
Furthermore, when $i=4$, $F_4$ collapses the non-simple part in the first and the second coordinates (the part that gets traversed back and forth along $\gamma$) of the loop $\bar{b}_2^{\scparallel}(t)$ down onto the lower left corner point $h$\/ of the left square in Figure~\ref{fig: disk}.  

Finally we check that the image of each $F_i$ is contained in the complement $Z''(1/q,m/r)\setminus\nu\bar{\Sigma}_2$.  
First we look at the case when $i=1$ or $2$.  Then either the first or the second coordinate of $F_i$ is 0.  From this, we deduce that the image of $F_i$ lies outside the $\epsilon$-tubular neighborhood of the $T_{\beta}$ part of $\bar{\Sigma}_2$ and also outside the 4-balls of radius $\rho$\/ centered at $x_2\times y_0$ and $x_3\times y_0$, wherein the resolution and the blow-up take place.  Since the third and the fourth coordinates of $F_i$ lie in the image of path $\Delta$, the image of $F_i$ lies outside the $\epsilon$-tubular neighborhood of the $T^2 \times y_0$ part of $\bar{\Sigma}_2$.  

We can argue similarly when $i=4$.  The first and the second coordinates of $F_4$ lie in the image of path $\gamma$, and hence the image of $F_4$ lies outside the $\epsilon$-tubular neighborhood of the $T_{\beta}$ part of $\bar{\Sigma}_2$ and outside the 4-balls of radius $\rho$\/ centered at $x_2\times y_0$ and $x_3\times y_0$.  
Since the third coordinate of $F_4$ lies in the interval $[0,\frac{1}{2}-\epsilon]$, the image of $F_4$ lies outside the $\epsilon$-tubular neighborhood of the $T^2 \times y_0$ part of $\bar{\Sigma}_2$.  
\end{proof}

\section{Fundamental group calculations}
\label{sec: pi_1}

The following theorem will serve as a cornerstone for many of the fundamental group calculations in later sections.  

\begin{thm}\label{thm:Z''}
There exists a nonnegative integer\/ $s$ such that\/ $\pi_1(Z''(1/q,m/r)\setminus\nu\bar{\Sigma}_2)$ is a quotient of the following finitely presented group
\begin{eqnarray}\label{pi_1(Sigma_2' complement in Z''(m))}
\langle
\alpha_1,\alpha_2,\alpha_3,\alpha_4, g_1, \dots, g_s 
&\mid& \alpha_3^q=[\alpha_1^{-1},\alpha_4^{-1}],\, \alpha_4^r=[\alpha_1,\alpha_3^{-1}]^m,\\
&&[\alpha_2,\alpha_3]=[\alpha_2,\alpha_4]=1
\rangle.  \nonumber
\end{eqnarray}
In\/ $\pi_1(Z''(1/q,m/r)\setminus\nu\bar{\Sigma}_2)$, 
the images of\/ $\alpha_1,\dots,\alpha_4$ are exactly the generators described at the beginning of Section\/~$\ref{sec: inclusion}$, and 
the images of\/ $g_1,\dots,g_s$ are elements of the subgroup normally generated by the image of\/ $[\alpha_3,\alpha_4]$.  
\end{thm}

\begin{proof}
Note that the immersed genus 2 surface that is the resolution of one intersection point between $T^2 \times y_0$ and $T_{\beta}$, i.e., the preimage of $\bar{\Sigma}_2$ before the blow-up, lies in the union $(T^2 \times D^2) \cup (D^2 \times T^2)$, where the $D^2$'s are the shaded disks in Figure~\ref{fig:squares}.  Thus it is easy to check that the tori 
\begin{eqnarray*}
\alpha_2\times\alpha_3 &=& \{(0,s,t,0)\mid 0\leq s,t \leq 1\}, \\
\alpha_2\times\alpha_4 &=& \{(0,s,0,t)\mid 0\leq s,t \leq 1\}
\end{eqnarray*}
are disjoint from $\bar{\Sigma}_2$ and the two surgery tori $\alpha_2' \times\alpha_3'$ and $\alpha_2''\times\alpha_4'$.  Hence we immediately obtain the commutator relations $[\alpha_2,\alpha_3]=[\alpha_2,\alpha_4]=1$.  As described at the beginning of Section~\ref{sec: inclusion}, each generator $\alpha_i$ is represented by a loop based at $h\times w_1$ which is the conjugate of a loop lying on the boundaries of two squares by the path $h\times\Delta$.  

The relation $\alpha_3^q=[\alpha_1^{-1},\alpha_4^{-1}]$ follows from the definition of $(\alpha_2'\times\alpha_3',\alpha_3',-1/q)$ surgery, once we realize that $\{\alpha_2,\alpha_3; [\alpha_1^{-1},\alpha_4^{-1}]\}$ is a ``distinguished triple of loops'' in the terminology of \cite{FPS}.
We can view $\alpha_3'$ as a loop based at $h\times w_1$ that is the composition of five paths: $h\times\Delta$, a bit of the path $t\mapsto (1-t,0,0,1-t)$ from $h\times k=(1,0,0,1)$ to $h'\times k'$, the original unbased loop $h'\times\alpha_3'$ shown in Figure~\ref{fig:squares}, a bit of the path $t\mapsto (t,0,0,t)$ from $h'\times k'$ to $h\times k$, and $h\times\Delta^{-1}$ in that order.  
As explained in \cite{FPS}, we can easily construct a based homotopy from $\alpha_3'$ to $\alpha_3$.  The image of such homotopy under the projection to the $\alpha_3\times\alpha_4$ torus factor, minus the $\Delta$ component, forms the thin rectangle whose sides lie on $\alpha_3$, $\alpha_3'$ and two copies of $\alpha_4$ in Figure~\ref{fig:squares}.  Similarly there is a rectangular homotopy between based paths $\alpha_2'$ and $\alpha_2$.  These homotopies can easily be made disjoint from $\nu\bar{\Sigma}_2$.  The torus $\alpha_1\times\alpha_4=\{(s,0,0,t)\mid 0\leq s,t \leq 1\}$ is disjoint from $\bar{\Sigma}_2$ and $\alpha_2''\times\alpha_4'$, but intersects the surgery torus $\alpha_2'\times\alpha_3'$ once at $h'\times k'$.  Hence the commutator $[\alpha_1^{-1},\alpha_4^{-1}]$ is a meridian of $\alpha_2'\times\alpha_3'$.  

The relation $\alpha_4^r=[\alpha_1,\alpha_3^{-1}]^m$ follows from the definition of $(\alpha_2''\times\alpha_4',\alpha_4',-m/r)$ surgery, once we check that 
the pair $\{\alpha_4;[\alpha_1,\alpha_3^{-1}]\}$ can be extended to a distinguished triple of loops.  As before, it is easy to construct a rectangular homotopy from $\alpha_4'$ to $\alpha_4$ whose image is disjoint from $\nu\bar{\Sigma}_2$.  

The torus $\alpha_1\times\alpha_3=\{(s,0,t,0)\mid 0\leq s,t \leq 1\}$ is disjoint from $\bar{\Sigma}_2$ and $\alpha_2'\times\alpha_3'$, but intersects the surgery torus $\alpha_2''\times\alpha_4'$ once.  The observations in \cite{BK:1-3, FPS} regarding orientation imply that the meridian of $\alpha_2''\times\alpha_4'$ in the distinguished triple should be either $[\alpha_1,\alpha_3^{-1}]$ or $[\alpha_1^{-1},\alpha_3^{-1}]$ if we include $\alpha_4$ in the distinguished triple.  As we traverse the path $t\mapsto (1-t,0)$ inside the left torus of Figure~\ref{fig:squares} in the negative $\alpha_1$-direction starting from $h$, we will have to cross $\alpha_2'$ before we reach $\alpha_2''$.  This obstruction implies that the path from $h \times k$\/ to the boundary 3-torus $\partial(\nu(\alpha_2''\times\alpha_4'))$ should be traversed in the positive $\alpha_1$-direction.  It follows that the correct choice of meridian in our distinguished triple should be $[\alpha_1,\alpha_3^{-1!
 }]$ (cf.\ \cite{BK:1-3, FPS}).  Thus we have shown that $\alpha_4^r=[\alpha_1,\alpha_3^{-1}]^m$ holds.  Since the surgery was along the $\alpha_4'$ curve, we do not actually need to find a concrete expression for the third member of our distinguished triple (the Lagrangian push-off of $\alpha_2''$) which will involve meridians of $\bar{\Sigma}_2$.  In fact, for all applications in this paper, it will be enough to know only that $\alpha_4^r=[\alpha_1^{\varepsilon_1},\alpha_3^{\varepsilon_3}]^m$ holds with $\varepsilon_1,\varepsilon_3\in\{\pm 1\}$.  

Finally, $[\alpha_1\times\alpha_2]\cdot[\bar{\Sigma}_2]=2$ and $[\alpha_3\times\alpha_4]\cdot[\bar{\Sigma}_2]=1$, so the relations $[\alpha_1,\alpha_2]=1$ and $[\alpha_3,\alpha_4]=1$ in $\pi_1(Z''(1/q,m/r))$ no longer hold in $\pi_1(Z''(1/q,m/r)\setminus\nu\bar{\Sigma}_2)$.   The commutator $[\alpha_3,\alpha_4]$ is a meridian of $\bar{\Sigma}_2$.  The blow-up process does not introduce any new generators in the fundamental group.  The extra generators $g_1,\dots,g_s$ correspond to meridians  of $\bar{\Sigma}_2$, and they do not appear in the relations of (\ref{pi_1(Sigma_2' complement in Z''(m))}). 
\end{proof}

\begin{cor}\label{cor:Z'}
There exists a nonnegative integer\/ $t$ such that $\pi_1(Z'\setminus\nu\bar{\Sigma}_2)$ is a quotient of the following finitely presented group
\begin{eqnarray}\label{pi_1(Sigma_2' complement)}
\langle
\alpha_1,\alpha_2,\alpha_3,\alpha_4, h_1, \dots, h_t 
&\mid& \alpha_3=[\alpha_1^{-1},\alpha_4^{-1}],\, [\alpha_1,\alpha_3]=1,\\
&& [\alpha_2,\alpha_3]=[\alpha_2,\alpha_4]=1
\rangle. \nonumber
\end{eqnarray}
In\/ $\pi_1(Z'\setminus\nu\bar{\Sigma}_2)$, the images of\/ $h_1,\dots,h_t$ are elements of the subgroup normally generated by the image of\/ $[\alpha_3,\alpha_4]$.  
There exists an embedding\/ $\bar{\Sigma}_2^{\scparallel}\hookrightarrow Z'\setminus\nu\bar{\Sigma}_2$ such that 
the induced homomorphism on\/ $\pi_1$ maps the generators\/ $\bar{a}_1^{\scparallel}$, $\bar{b}_1^{\scparallel}$ and\/ $\bar{b}_2^{\scparallel}$ of\/ $\pi_1(\bar{\Sigma}_2^{\scparallel})$ to the images of\/ $\alpha_1$, $\alpha_2$ and\/ $\alpha_4$, respectively.  
\end{cor}

\begin{proof}
Note that $Z''(1/1,1/0)=Z'$ and hence we can set $q=m=1$ and $r=0$ in Lemma~\ref{lem: parallel} and Theorem~\ref{thm:Z''}.  Clearly $[\alpha_1,\alpha_3^{-1}]=1$ implies $[\alpha_1,\alpha_3]=1$.  Geometrically, the torus $\alpha_1\times\alpha_3=\{(s,0,t,0)\mid 0\leq s,t \leq 1\}$ is disjoint from both $\alpha_2'\times\alpha_3'$ and $\bar{\Sigma}_2$, so it survives into $Z'\setminus\nu\bar{\Sigma}_2$.   
\end{proof}

\section{Construction of exotic $\CP\#2\CPb$} 
\label{sec:exotic 1-2}

For each pair of integers $p \geq 0$ and $q\geq 0$, let $Y_1(1/p,1/q)$ be the result of the following 4 torus surgeries on $\Sigma_2\times T^2$:
\begin{equation}\label{eq: Luttinger surgeries for Y_1(m)} 
(a_1' \times c', a_1', -1), \ \ (b_1' \times c'', b_1', -1),\ \
(a_2' \times c', c', +1/p), \ \ (a_2'' \times d', d', +1/q).
\end{equation}
Here, $a_i,b_i$ ($i=1,2$) and $c,d$\/ denote the standard generators of $\pi_1(\Sigma_2)$ and $\pi_1(T^2)$, respectively.  The surgeries in (\ref{eq: Luttinger surgeries for Y_1(m)}) are all Luttinger surgeries when $p\geq 1$ and $q\geq 1$.  Hence $Y_1(1/p,1/q)$ is a minimal symplectic 4-manifold for all $p\geq 0$ and $q\geq 0$.  

$\pi_1(Y_1(1/p,1/q))$ is generated by 
$a_i,b_i$ ($i=1,2$) and $c,d$.  The following relations hold in 
$\pi_1(Y_1(1/p,1/q))$:
\begin{gather}\label{Luttinger relations for Y_1(m)}
[b_1^{-1},d^{-1}]=a_1,\ \  [a_1^{-1},d]=b_1,\ \
[d^{-1},b_2^{-1}]=c^p,\ \ [c^{-1},b_2]=d^q,\\ \nonumber
[a_1,c]=1,\ \  [b_1,c]=1,\ \ [a_2,c]=1,\ \  [a_2,d]=1,\\ \nonumber
[a_1,b_1][a_2,b_2]=1,\ \ [c,d]=1.
\end{gather}

Let $\Sigma_2=\Sigma_2\times(\frac{1}{2},\frac{1}{2})\subset Y_1(1/p,1/q)$ be a genus 2 symplectic submanifold.  Now
let $p=q=r=1$ and take the normal connected sum 
\begin{equation*}
X_1(m)=Y_{1}(1,1)\#_{\psi}Z''(1,m)
\end{equation*}
using an orientation reversing diffeomorphism $\psi:\partial(\nu\Sigma_{2})\rightarrow \partial(\nu\bar{\Sigma}_2)$.  Here, $\nu$ denotes tubular neighborhoods.  
$\psi$\/ restricts to orientation preserving diffeomorphism on parallel genus 2 surfaces and complex conjugation on the meridian circles.

As in the proof of Theorem~\ref{thm:Z''}, we can choose the base point $z_1$ of $\pi_1(Y_1(1,1))$ on $\partial(\nu\Sigma_2)$ such that $\pi_1(Y_1(1,1)\setminus\nu\Sigma_2,z_1)$ is normally generated by $a_i,b_i$ ($i=1,2$) and $c,d$.  These generators are now represented by loops based at $z_1$ that are $\Delta_1$-conjugates of standard loops based at the product $h_1\times k_1$ of corner points of an octagon and a square, where $\Delta_1$ is a fixed path from $z_1$ to $h_1\times k_1$.  Since $\Sigma_2=\Sigma_2\times(\frac{1}{2},\frac{1}{2})$ is disjoint from the neighborhoods of four Luttinger surgery tori, 
all the relations in (\ref{Luttinger relations for Y_1(m)}) continue to hold in $\pi_1(Y_1(1,1)\setminus\nu\Sigma_2)$ except for
$[c,d]=1$.  This commutator is no longer trivial and now represents a meridian of $\Sigma_{2}$ in $\pi_1(Y_1(1,1)\setminus\nu\Sigma_2)$.  

We now require $\psi_{\ast}$ to map the generators of $\pi_1$ as follows:
\begin{equation}\label{eq:psi mapsto for X_1(m)}
a_i \mapsto \bar{a}_i^{\scparallel}, \ \ 
b_i \mapsto \bar{b}_i^{\scparallel}, \ \ 
i = 1,2.  
\end{equation}
Note that $X_1(m)$ is symplectic when $m=1$ (cf.\ \cite{gompf}). 

\begin{lem}\label{lem:1-2}
The set\/ $S_{1,2}=\{X_1(m)\mid m\geq 1\}$ consists of irreducible\/ $4$-manifolds that are homeomorphic to $\CP\# 2\CPb$.  Moreover, $S_{1,2}$ contains an infinite subset consisting of pairwise non-diffeomorphic non-symplectic\/ $4$-manifolds.
\end{lem}

\begin{proof}
We have 
\begin{eqnarray*}
e(X_1(m)) &=& e(Y_{1}(1,1)) + e(Z''(1,m)) -2 e(\Sigma_{2}) = 0 + 1 + 4 = 5,\\
\sigma(X_1(m)) &=& \sigma (Y_{1}(1,1)) + \sigma(Z''(1,m)) = 0 + (-1) = -1.
\end{eqnarray*}
>From Freedman's theorem (cf.\ \cite{freedman}), we conclude that $X_1(m)$ is homeomorphic to $\CP\# 2\CPb$, once we show that $\pi_1(X_1(m))=1$.   From Seifert-Van Kampen theorem,
we can deduce that $\pi_1(X_1(m))$ is a quotient of the following group:  
\begin{equation}\label{pi_1(X_1(m))}
\frac{\pi_1(Y_{1}(1,1)\setminus\nu\Sigma_{2})\ast
\pi_1(Z''(1,m)\setminus\nu\bar{\Sigma}_2)}{\langle a_1=\alpha_1,\,
b_1=\alpha_2,\, 
b_2=\alpha_4,\, 
\mu(\Sigma_{2})=\mu(\bar{\Sigma}_2)^{-1} \rangle},
\end{equation}

In (\ref{pi_1(X_1(m))}), we have $a_1=[b_1^{-1},d^{-1}]=[b_1^{-1},[c^{-1},b_2]^{-1}]=[b_1^{-1},[b_2,c^{-1}]]$.  
Since $[\alpha_2,\alpha_4]=1$ in (\ref{pi_1(Sigma_2' complement in Z''(m))}), we conclude that $[b_1, b_2]=1$.  Since we also have $[b_1,c]=1$ in (\ref{Luttinger relations for Y_1(m)}), we easily deduce that $a_1=1$.  
>From $a_1=1$ and (\ref{Luttinger relations for Y_1(m)}), we conclude that $b_1=\alpha_2=1$.  Plugging $\alpha_1=a_1=1$ into (\ref{pi_1(Sigma_2' complement in Z''(m))}), we obtain $\alpha_3=\alpha_4 =1$.  This in turn implies that $[\alpha_3,\alpha_4]$, a meridian of $\bar{\Sigma}_2$, is trivial and hence the generators $g_1,\dots,g_s$ coming from (\ref{pi_1(Sigma_2' complement in Z''(m))}) are all trivial as well.  Next, $\alpha_4=1$ implies that $b_2=1$, and then it follows from (\ref{Luttinger relations for Y_1(m)}) that $c=d=1$.  Finally, since we now know that the meridians of $\bar{\Sigma}_2$ are all trivial, (\ref{eq: embedding of Sigma_2'}) implies that $a_2=\alpha_3^2=1$.  In conclusion, (\ref{pi_1(X_1(m))}) is the trivial group and thus $\pi_1(X_1(m))=1$.   

To prove irreducibility, we first observe that $Y_1(1,1)$ is minimal and that the only $-1$ sphere in $Z''(1,m)$ is the exceptional sphere $E$\/ of the blow-up by Corollary~3 in \cite{tjli}.  Note that $E$\/ intersects $\bar{\Sigma}_2$ twice in $Z''(1,m)$ and thus there is no $-1$ sphere in $Z''(1,m)\setminus\nu\bar{\Sigma}_2$.  
It follows from Usher's theorem in \cite{usher} that the symplectic normal connected sum 
$X_1(1)$ is symplectically minimal.  
Recall that symplectic minimality implies irreducibility when $\pi_1=1$ by 
the theorem of Hamilton and Kotschick in \cite{HK}.   Hence $X_1(1)$ is irreducible.  
Since $X_1(m)$ is obtained from $X_1(1)$ by performing a  $1/(m-1)$ surgery on a  nullhomologous torus (cf.\ \cite{ABP, FPS}), $X_1(m)$ must be irreducible as well.  

To prove that infinitely many $X_1(m)$'s are pairwise non-diffeomorphic, we need to view $X_1(m)$ as the result of $5$ Luttinger surgeries and a single $m$\/ torus surgery on the symplectic normal connected sum $(\Sigma_2 \times T^2)\#_{\psi}(T^4\#\CPb)$, which is a minimal symplectic $4$-manifold with $b_2^+>1$.  In other words, we view 
$X_1(m)$ as being ``reverse-engineered'' from $(\Sigma_2 \times T^2)\#_{\psi}(T^4\#\CPb)$ in the spirit of \cite{FPS, stern}.  
We can then compute the Seiberg-Witten invariants of $X_1(m)$ and check that infinitely many of them are distinct by applying exactly the same argument as in \cite{ABP, FPS} using the product formulas in \cite{MMS}.   

Moreover, we can also observe that the values of the Seiberg-Witten invariants of $X_1(m)$ grow arbitrarily large as $m\rightarrow\infty$, and in particular these values cannot be $\pm 1$ when $m$\/ is large.  
Since the value of the Seiberg-Witten invariant on the canonical class of a symplectic $4$-manifold is always $\pm 1$ by Taubes's theorem (cf.\ \cite{taubes}), $X_1(m)$ cannot be symplectic when $m$\/ is large.  
\end{proof}

For the remainder of the section, we look at some interesting constructions of $4$-manifolds with cyclic fundamental groups.  
Consider the following new family of normal connected sums 
$$\tilde{X}_1(p,m)=Y_{1}(1/p,1)\#_{\psi}Z''(1,m),$$
where the gluing diffeomorphism $\psi: \partial(\nu\Sigma_{2})\rightarrow \partial(\nu\bar{\Sigma}_2)$ still satisfies (\ref{eq:psi mapsto for X_1(m)}).  
Note that $\tilde{X}_1(p,m)$ is symplectic when $m=1$.  

The proof of Lemma~\ref{lem:1-2} goes through almost verbatim for $\tilde{X}_1(p,m)$.   In particular, $\tilde{X}_1(p,m)$ has the same Euler characteristic and signature as $\CP\#2\CPb$.  The only difference is that the generator $c$\/ is no longer trivial when $p\neq 1$.  It is not hard to see that the order of $c$\/ is unaffected by the relations coming from $\pi_1(Z''(1,m)\setminus\nu\bar{\Sigma}_2)$.  In summary, we have 
\begin{equation*}
\pi_1(\tilde{X}_1(p,m))=\left\{\begin{array}{cll}
\Z & \textrm{ if } &  p=0, \\
\Z/p & \textrm{ if } & p\geq 2. 
\end{array}\right.
\end{equation*}

\begin{rmk}
Note that the Euler characteristic of each irreducible symplectic\/ $4$-manifold $\tilde{X}_1(p,1)$ is $5$, which is one less than the smallest examples constructed in \cite{BK:other constructions} with fundamental group equal to $\Z$ or $\Z/p$.  
The authors will deal with non-cyclic fundamental groups in a forthcoming paper.  
\end{rmk}

\section{Construction of exotic $(2n-1)\CP\# 2n\CPb$ for $n\geq 2$}
\label{sec:exotic 3-4}

Let $Y_{n}(m)$ be the manifold from Section~\ref{sec:(2n-3)(S^2 x S^2)} with 
a genus 2 surface $\Sigma_2$ sitting inside it.  We take the normal connected sum $X_n(m)=Y_{n}(m)\#_{\psi}Z'$ using a diffeomorphism  $\psi:\partial(\nu\Sigma_2)\to \partial(\nu\bar{\Sigma}_2)$.  We require $\psi_{\ast}$ to map the generators of $\pi_1$ as follows:
\begin{equation}\label{eq:psi mapsto}
a_i \mapsto \bar{a}_i^{\scparallel}, \ \ 
b_i \mapsto \bar{b}_i^{\scparallel}, \ \ 
i = 1,2.  
\end{equation}
Note that $X_n(m)$ is symplectic when $m=1$. 

\begin{lem}\label{lem:(2n-1)-2n}
For each integer\/ $n\geq 2$, the set\/ $S_{2n-1,2n}=\{X_n(m)\mid m\geq 1\}$ consists of irreducible\/ $4$-manifolds that are homeomorphic to $(2n-1)\CP\# 2n\CPb$.
Moreover, $S_{2n-1,2n}$ contains an infinite subset consisting of pairwise non-diffeomorphic non-symplectic\/ $4$-manifolds.
\end{lem}

\begin{proof}
We can easily compute that
\begin{eqnarray*}
e(X_n(m)) &=& e(Y_{n}(m)) + e(Z') -2 e(\Sigma_{2}) = (4n-4) + 1 + 4 = 4n+1,\\
\sigma(X_n(m)) &=& \sigma (Y_{n}(m)) + \sigma(Z') = 0 + (-1) = -1.
\end{eqnarray*}
>From Freedman's theorem (cf.\ \cite{freedman}), we conclude that $X_n(m)$ is homeomorphic to $(2n-1)\CP\# 2n\CPb$, once we show that $\pi_1(X_n(m))=1$.   From Seifert-Van Kampen theorem,
we can deduce that $\pi_1(X_n(m))$ is a quotient of the following group:  
\begin{equation}\label{pi_1(X_n(m))}
\frac{\pi_1(Y_{n}(m)\setminus\nu\Sigma_{2})\ast
\pi_1(Z'\setminus\nu\bar{\Sigma}_2)}{\langle a_1=\alpha_1,\,
b_1=\alpha_2,\, 
b_2=\alpha_4,\, 
\mu(\Sigma_{2})=\mu(\bar{\Sigma}_2)^{-1} \rangle}.
\end{equation}
All the relations in $\pi_1(Y_{n}(m))$ listed in (\ref{Luttinger relations}) of Section~\ref{sec:(2n-3)(S^2 x S^2)} continue to hold in (\ref{pi_1(X_n(m))}) except possibly for
$\prod_{j=1}^n[c_j,d_j]=1$.  This product may no longer be trivial and now represents a meridian of $\Sigma_{2}$.  

In (\ref{pi_1(X_n(m))}), we have $a_1=[b_1^{-1},d_1^{-1}]=[b_1^{-1},[c_1^{-1},b_2]^{-1}]=[b_1^{-1},[b_2,c_1^{-1}]]$.  
Since $[\alpha_2,\alpha_4]=1$ in (\ref{pi_1(Sigma_2' complement)}), we conclude that $[b_1,b_2]=1$.  Since we also have $[b_1,c_1]=1$ from (\ref{Luttinger relations}), we easily deduce that $a_1=1$.  
>From $a_1=1$, we can rapidly kill all the generators of the form $a_i,\ b_i,\ c_j,\ d_j$.  Consequently, when we pass from (\ref{pi_1(Sigma_2' complement)}) to $\pi_1(X_n(m))$, all the generators of the form $\alpha_i$ also die.  For example, $\alpha_3=[\alpha_1^{-1},\alpha_4^{-1}]=[a_1^{-1},b_2^{-1}]=1$.  Finally note that $[\alpha_3,\alpha_4]$, a meridian of $\bar{\Sigma}_2$, is trivial and hence the generators $h_1,\dots,h_t$ die as well.  In conclusion, (\ref{pi_1(X_n(m))}) is the trivial group and consequently $\pi_1(X_n(m))=1$.  The same arguments as in the proof of Lemma~\ref{lem:1-2} show that $X_n(m)$ are all irreducible and infinitely many of them are non-symplectic and pairwise non-diffeomorphic.  
\end{proof}

\begin{rmk}
By taking the normal connected sum $Y_n(m)\#_{\psi}(T^4\#\CPb)$ using the diffeomorphism 
(\ref{eq:psi mapsto}), i.e., without performing the Luttinger surgery on $T^4\#\CPb$, we obtain an infinite family of pairwise homeomorphic (cf.\ \cite{hambleton-kreck}) 4-manifolds that have $\pi_1=\Z/2$ and the same rational cohomology ring as $(2n-1)\CP\# 2n \CPb$ for $n\geq 2$. 
\end{rmk}

\section{Construction of exotic $\CP\#4\CPb$ and $\CP\#6\CPb$}
\label{sec:exotic 1-4}

To construct exotic $\CP\#4\CPb$, we start with a genus 2 symplectic surface of self-intersection 0 in $T^4\#2\CPb$.   Assume that $T^4=T^2 \times T^2$ is equipped with a product symplectic form.  Let $\alpha_i$ be as in Section~\ref{sec:braided surface}.  
Symplectically resolve the intersection between two symplectic tori $\alpha_1\times\alpha_2\times (\frac{1}{2},\frac{1}{2})$ and $(\frac{1}{2},\frac{1}{2})\times\alpha_3\times\alpha_4$ in $T^4$ to obtain a genus 2 surface of self-intersection 2 in $T^4$.  Now symplectically blow up twice to obtain a genus 2 symplectic surface $\hat{\Sigma}_2$ of self-intersection 0 in $T^4\#2\CPb$.   Given a pair of positive integers $q$\/ and $r$, let $M(1/q,1/r)$ be the result of the following two Luttinger surgeries on $T^4\#2\CPb$:
\begin{equation}\label{eq:  Luttinger surgeries for M(1/q,1/r)}
(\alpha_1' \times \alpha_3', \alpha_1', -1/q), \ \ 
(\alpha_2' \times \alpha_3'', \alpha_2', -1/r).
\end{equation} 
$M(1/q,1/r)$ is diffeomorphic to $Z''(1/q,1/r)\#\CPb$.  
Note that $\hat{\Sigma}_2$ lies away from the Lagrangian tori used in the Luttinger surgeries in (\ref{eq: Luttinger surgeries for M(1/q,1/r)}) and hence $\hat{\Sigma}_2$ descends to a symplectic submanifold in $M(1/q,1/r)$.   We may assume that the inclusion $\hat{\Sigma}_2\hookrightarrow M(1/q,1/r)$ maps the standard generators of $\pi_1(\hat{\Sigma}_2)$ as follows:
\begin{equation*}
\hat{a}_1 \mapsto \alpha_1, \ \ 
\hat{b}_1 \mapsto \alpha_2, \ \
\hat{a}_2 \mapsto \alpha_3, \ \ 
\hat{b}_2 \mapsto \alpha_4.  
\end{equation*}

\begin{lem}
With\/ $\hat{\Sigma}_2\subset M(1/q,1/r)$ as above, we have 
\begin{eqnarray}
\label{eq: pi_1(Sigma_2'' complement)}
\hspace{15pt}
\pi_1(M(1/q,1/r)\setminus\nu\hat{\Sigma}_2) &=& 
\langle
\alpha_1,\alpha_2,\alpha_3,\alpha_4 \mid
\alpha_1^q=[\alpha_2^{-1},\alpha_4^{-1}],\\ \nonumber
&&\hspace{4pt} \alpha_2^r=[\alpha_1^{-1},\alpha_4],\,
[\alpha_1,\alpha_3]=[\alpha_2,\alpha_3]=1,\, \\ \nonumber
&&\hspace{4pt} [\alpha_1,\alpha_2]=[\alpha_3,\alpha_4]=1
\rangle.
\end{eqnarray}
\end{lem}

\begin{proof}
Note that $[\alpha_1\times\alpha_2]\cdot[\hat{\Sigma}_2]=[\alpha_3\times\alpha_4]\cdot[\hat{\Sigma}_2]=1$.  The meridian of $\hat{\Sigma}_2$ is trivial in $\pi_1(M(1/q,1/r)\setminus\nu\hat{\Sigma}_2)$ since it bounds a punctured exceptional sphere of a blow-up.  Hence we must have $[\alpha_1,\alpha_2]=[\alpha_3,\alpha_4]=1$.  The other relations coming from Luttinger surgeries can be derived in exactly the same way as before.  
\end{proof}

Let $q=r=1$ and take the normal connected sum $V(m)=M(1,1)\#_{\psi}Z''(1,m)$, where the gluing diffeomorphism $\psi:\partial(\nu\hat{\Sigma}_2)\rightarrow\partial(\nu\bar{\Sigma}_2)$ maps the generators of $\pi_1(\hat{\Sigma}_2^{\scparallel})$ as follows:
\begin{equation}\label{eq:psi mapsto for 1-4}
\hat{a}_i^{\scparallel} \mapsto \bar{a}_i^{\scparallel}, \ \ 
\hat{b}_i^{\scparallel} \mapsto \bar{b}_i^{\scparallel}, \ \
i=1,2.
\end{equation}
Note that $V(m)$ is symplectic when $m=1$.  

\begin{lem}\label{lem:V(m)}
The set\/ $S_{1,4}=\{ V(m)\mid m\geq 1\}$ consists of irreducible\/ $4$-manifolds that are homeomorphic to $\CP\#4\CPb$.  Moreover, $S_{1,4}$ contains an infinite subset consisting of pairwise non-diffeomorphic non-symplectic\/ $4$-manifolds.  
\end{lem}

\begin{proof}
We compute that 
\begin{eqnarray*}
e(V(m)) &=& e(M(1,1)) + e(Z''(1,m)) -2 e(\hat{\Sigma}_{2}) = 2 + 1 + 4 = 7,\\
\sigma(V(m)) &=& \sigma (M(1,1)) + \sigma(Z''(1,m)) = (-2) + (-1) = -3.
\end{eqnarray*}
Freedman's theorem (cf.\ \cite{freedman}) implies that $V(m)$ is homeomorphic to $\CP\# 4\CPb$, once we show that $\pi_1(V(m))=1$.   From Seifert-Van Kampen theorem, we deduce that $\pi_1(V(m))$ is a quotient of the following finitely presented group:  
\begin{eqnarray}\label{eq: pi_1(V(m))}
\hspace{15pt}
\langle
\alpha_1,\alpha_2,\alpha_3,\alpha_4, g_1, \dots, g_s 
&\mid& \alpha_3=[\alpha_1^{-1},\alpha_4^{-1}],\, \alpha_4=[\alpha_1,\alpha_3^{-1}]^m,\\ \nonumber
&&[\alpha_2,\alpha_3]=[\alpha_2,\alpha_4]=1,\\ \nonumber
&&\alpha_1=[\alpha_2^{-1},\alpha_4^{-1}],\, \alpha_2=[\alpha_1^{-1},\alpha_4]
\rangle,
\end{eqnarray}
where $g_1, \dots, g_s$ are all conjugate to $[\alpha_3,\alpha_4]$ in $\pi_1(V(m))$.  
It is easy to see that $\alpha_1=1$ in (\ref{eq: pi_1(V(m))}), which implies that all other generators of (\ref{eq: pi_1(V(m))}) are also trivial in $\pi_1(V(m))$.  It follows that $\pi_1(V(m))=1$.   The same arguments as in the proof of Lemma~\ref{lem:1-2} show that $V(m)$ are all irreducible and infinitely many of them are non-symplectic and pairwise non-diffeomorphic.  
\end{proof}

To construct exotic $\CP\#6\CPb$, we proceed as follows.  First recall that there is a genus $2$ symplectic surface $\widetilde{\Sigma}_2$ of self-intersection $0$ in $(T^2\times S^2)\#4\CPb$.  $\widetilde{\Sigma}_2$ is obtained by symplectically resolving the intersections between $\{{\rm pt}\}\times S^2$ and two parallel copies of $T^2\times\{{\rm pt}\}$ and then symplectically blowing up $4$ times.  Let us denote the standard generators of $\pi_1(\widetilde{\Sigma}_2)$ and $\pi_1((T^2\times S^2)\#4\CPb)\cong \pi_1(T^2)$ by $\widetilde{a}_i,\widetilde{b}_i$ ($i=1,2$) and $c,d,$ respectively.  We can assume that the inclusion $\widetilde{\Sigma}_2 \hookrightarrow (T^2\times S^2)\#4\CPb$ maps the generators as follows:  
\begin{equation}\label{eq: embedding of tilde{Sigma}_2}
\widetilde{a}_1 \mapsto c, \ \ 
\widetilde{b}_1 \mapsto d, \ \
\widetilde{a}_2 \mapsto c^{-1}, \ \ 
\widetilde{b}_2 \mapsto d^{-1}.  
\end{equation}
It is easy to see that $\pi_1(((T^2\times S^2)\#4\CPb)\setminus\nu\widetilde{\Sigma}_2)\cong \pi_1((T^2\times S^2)\#4\CPb)\cong \Z^2$ since each exceptional sphere intersects $\widetilde{\Sigma}_2$ once and hence the meridian of $\widetilde{\Sigma}_2$ is nullhomotopic in the complement of $\widetilde{\Sigma}_2$.  

Now take the normal connected sum $W(m)=((T^2\times S^2)\#4\CPb)\#_{\psi}Z''(1,m)$, where the gluing diffeomorphism $\psi: \partial(\nu \widetilde{\Sigma}_2) \rightarrow \partial(\nu \bar{\Sigma}_2)$ maps
the generators of $\pi_1(\widetilde{\Sigma}_2^{\scparallel})$ as follows:  
\begin{equation}\label{eq:psi mapsto for W(m)}
\widetilde{a}_i^{\scparallel} \mapsto \bar{a}_i^{\scparallel}, \ \ 
\widetilde{b}_i^{\scparallel} \mapsto \bar{b}_i^{\scparallel}, \ \ 
i = 1,2.  
\end{equation}
Note that $W(m)$ is symplectic when $m=1$.  

\begin{lem}\label{lem:W(m)}
The set\/ $S_{1,6}=\{W(m)\mid m\geq 1\}$ consists of irreducible\/ $4$-manifolds that are homeomorphic to $\CP\#6\CPb$.  Moreover, $S_{1,6}$ contains an infinite subset consisting of pairwise non-diffeomorphic non-symplectic\/ $4$-manifolds.
\end{lem}

\begin{proof}
We compute that 
\begin{eqnarray*}
e(W(m)) &=& e((T^2\times S^2)\#4\CPb) + e(Z''(1,m)) -2 e(\widetilde{\Sigma}_2) = 4 + 1 + 4 = 9,\\
\sigma(W(m)) &=& \sigma ((T^2\times S^2)\#4\CPb) + \sigma(Z''(1,m)) = (-4) + (-1) = -5.
\end{eqnarray*}
Freedman's theorem (cf.\ \cite{freedman}) implies that $W(m)$ is homeomorphic to $\CP\# 6\CPb$, once we show that $\pi_1(W(m))=1$.   From Seifert-Van Kampen theorem, we deduce that $\pi_1(W(m))$ is a quotient of the following finitely presented group:  
\begin{eqnarray}\label{eq: pi_1(W(m))}
\hspace{15pt}
\langle
\alpha_1,\alpha_2,\alpha_3,\alpha_4, g_1, \dots, g_s 
&\mid& \alpha_3=[\alpha_1^{-1},\alpha_4^{-1}],\, \alpha_4=[\alpha_1,\alpha_3^{-1}]^m, \\ \nonumber
&&[\alpha_2,\alpha_3]=[\alpha_2,\alpha_4]=1,\\ \nonumber
&&  \alpha_4 = \alpha_2^{-1},\,
[\alpha_1,\alpha_2] = 1
\rangle,
\end{eqnarray}
where $g_1, \dots, g_s$ are all conjugate to $[\alpha_3,\alpha_4]$ in $\pi_1(W(m))$.  
Note that we have $\alpha_3=[\alpha_1^{-1},\alpha_2]=1$ in (\ref{eq: pi_1(W(m))}), which easily implies that all other generators of (\ref{eq: pi_1(W(m))}) are also trivial in $\pi_1(W(m))$.  For example, since the meridian $[\alpha_3,\alpha_4]$ of $\bar{\Sigma}_2$ is trivial, (\ref{eq: embedding of Sigma_2'}), (\ref{eq: embedding of tilde{Sigma}_2}) and (\ref{eq:psi mapsto for W(m)}) imply that $\alpha_1=\alpha_3^{-2}=1$.  Hence $\pi_1(W(m))=1$.  
The same arguments as in the proof of Lemma~\ref{lem:1-2} show that $W(m)$ are all irreducible and infinitely many of them are non-symplectic and pairwise non-diffeomorphic.  
\end{proof}

An irreducible symplectic $4$-manifold homeomorphic to $\CP\#6\CPb$ was first constructed in \cite{SS}.  An infinite family of non-symplectic irreducible $4$-manifolds homeomorphic to $\CP\#6\CPb$ was first constructed in \cite{FS: double node}.   
Our exotic symplectic $4$-manifold $W(1)$ is new in the sense that it is known to contain a genus 2 symplectic surface of self-intersection 0.   This fact allows us to build an irreducible symplectic $4$-manifold homeomorphic to $3\CP\#10\CPb$ in Section~\ref{sec:exotic 3-k}.

\section{Construction of exotic $3\CP\#k\CPb$ for $k=6,8,10$}
\label{sec:exotic 3-k}

\begin{lem}\label{lem:3-k}
For $k=4,6,8,10$, there exist an irreducible symplectic\/ $4$-manifold and an infinite family of pairwise non-diffeomorphic irreducible non-symplectic\/ $4$-manifolds that are all homeomorphic to\/ $3\CP\#k\CPb$.
\end{lem}

\begin{proof}
Note that such an infinite family of exotic $3\CP\#4\CPb$'s was already constructed in Section~\ref{sec:exotic 3-4}.   The existence of exotic $3\CP\#k\CPb$ for $k=6,8,10$ follows at once from the existence of exotic $\CP\#(k -4)\CPb$ containing a square 0 genus 2 surface in Sections \ref{sec:exotic 1-2} and \ref{sec:exotic 1-4} by normally connect-summing with $T^4\#2\CPb$ along the proper transform $\hat{\Sigma}_2$ of the genus 2 resolution of $(\alpha_1\times\alpha_2) \cup (\alpha_3\times\alpha_4)$ (cf.\ Theorem 6 of \cite{ABP}).  Alternatively, we can normally connect-sum $Z'$ of Section~\ref{sec: Luttinger surgery} to exotic $\CP\#(k-3)\CPb$'s constructed in \cite{akhmedov,ABP,AP,BK:1-3,BK:other constructions} along genus 2 surfaces.  The details can be filled in as in the previous proofs and are left to the reader.  
\end{proof}

\begin{rmk}
Infinitely many non-symplectic irreducible $4$-manifolds homeomorphic to 
$3\CP\#8\CPb$ and $3\CP\#10\CPb$ were already constructed in \cite{jpark:3-8} and \cite{bdpark:3-n}, respectively.   The constructions in this section give the first known irreducible {\em symplectic}\/ $4$-manifolds homeomorphic to these.  
\end{rmk}

\section{Proof of Theorem~\ref{thm:wedge}}
\label{sec:proof of wedge thm}

For each pair of nonnegative integers $(\chi,c)$ satisfying (\ref{eq: chi-c inequality}), we will construct an odd minimal symplectic $4$-manifold $N$\/ containing a symplectic torus $T'$ of self-intersection $0$ such that 
\begin{equation*}
\chi_h(N) = \chi, \quad c_1^2(N) = c, 
\end{equation*}
and the inclusion induced homomorphism $\pi_1(T')\rightarrow \pi_1(N)$ is surjective.  We will also require that the meridian of $T'$ is trivial in $\pi_1(N\setminus\nu T')$.  This implies that\/ $\pi_1(N\setminus\nu T')\cong\pi_1(N)$.

Given such a pair $(N,T')$, we can define $Y$\/ to be the symplectic normal connected sum $X\#_{\psi}N$, where $\psi:\partial(\nu T)\rightarrow \partial(\nu T')$ denotes a suitable orientation reversing diffeomorphism.  Identities in (\ref{eq: chi_h and c_1^2}) are immediately satisfied.  $Y$\/ is odd since $N$\/ is odd.  If $X$\/ is minimal, then the minimality of $Y$\/ follows at once from Usher's theorem (cf.\ \cite{usher}).  We can compute $\pi_1(Y)$ by Seifert-Van Kampen theorem.  First note that the meridian of $T$\/ is trivial in $\pi_1(Y)$ since it gets identified with the inverse of the meridian of $T'$, which is trivial by our assumption in the previous paragraph.  Let $T^{\scparallel}$ and $(T')^{\scparallel}$ denote parallel copies of $T$\/ and $T'$ in the boundaries of tubular neighborhoods $\nu T\subset X$\/ and $\nu T'\subset N$, respectively.  The inclusion induced homomorphism $\pi_1(N\setminus\nu T')\rightarrow \pi_1(Y)$ is trivial since the generators of $!
 \pi_1((T')^{\scparallel})$ get identified with the generators of $\pi_1(T^{\scparallel})$, which are trivial in $\pi_1(Y)$ by the hypothesis of Theorem~\ref{thm:wedge}.  It follows that $\pi_1(Y) \cong \pi_1(X\setminus\nu T)/\langle \mu(T)\rangle \cong \pi_1(X)$.  

It will be enough to construct $(N,T')$ corresponding to the pairs $(\chi,c)$ equal to
\begin{equation}\label{eq: pairs to check}
(1,5),\  (2,9),\  (2,11),\  (2,13),\  \text{or}\  (\chi,8\chi-1)\  \text{with} \ \chi\geq 1.  
\end{equation}
Recall that the other remaining pairs $(\chi,c)$ have already been dealt with in \cite{ABBKP}.  Note that the pairs in (\ref{eq: pairs to check}) correspond to $(\chi_h(M),c_1^2(M))$, where $M$\/ is $\CP\#4\CPb$, $3\CP\#10\CPb$, $3\CP\#8\CPb$, $3\CP\#6\CP$, or $(2\chi-1)\CP\#2\chi\CPb$ with $\chi\geq 1$, respectively.  

To construct $(N,T')$ corresponding to each of these $(\chi,c)=(\chi_h(M),c_1^2(M))$, we proceed as follows.  In Sections~\ref{sec:exotic 1-2}--\ref{sec:exotic 3-k}, we have constructed a minimal symplectic $4$-manifold that is homeomorphic to $M$.  Each such minimal symplectic $4$-manifold can be viewed as the result of performing multiple Luttinger surgeries on a symplectic normal connected sum.  If we choose to forgo a single Luttinger surgery along a Lagrangian torus $T'$, then we will be left with a minimal symplectic $4$-manifold $N$\/ that has the same $\chi_h$ and $c_1^2$\/ as $M$, but such $N$\/ will no longer be simply-connected.  It is possible to choose a suitable Lagrangian torus $T'$ such that the meridian of $T'$ is trivial in $\pi_1(N\setminus\nu T')$ and the inclusion induced homomorphism $\pi_1(T')\rightarrow\pi_1(N)$ is surjective.  Table~\ref{table: forgo} lists a possible choice for each $(\chi,c)$ in (\ref{eq: pairs to check}).

\begin{table}[ht]
\caption{Possible choice of a Luttinger surgery to skip}
\begin{center}
\begin{tabular}{|c|c|}
\hline
$(\chi, c)$ & Luttinger surgery \\ \hline
$(1,7)$ & $(a_2' \times c', c', +1)$ in $Y_1(1,1)$ summand of $X_1(1)$ \\ \hline
$(\chi,8\chi-1)$ with $\chi\geq 2$ & $(a_2' \times c_1', c_1', +1)$ in $Y_{\chi}(1)$ summand of $X_{\chi}(1)$ \\ \hline
$(1,5)$ & $(\alpha_2'' \times\alpha_4', \alpha_4', -1)$ in $Z''(1,1)$ summand of $V(1)$ \\ \hline
$(2,13)$ & $(a_2' \times c', c', +1)$ in $X_1(1)$ summand \\ \hline
$(2,11)$ & $(\alpha_2'' \times\alpha_4', \alpha_4', -1)$ in $V(1)$ summand \\ \hline
$(2,9)$ & $(\alpha_2'' \times\alpha_4', \alpha_4', -1)$ in $W(1)$ summand \\
\hline 
\end{tabular}
\end{center}
\label{table: forgo}
\end{table}

For example, if we choose not to perform $(a_2' \times c', c', +1)$ Luttinger surgery in the $Y_1(1,1)$ summand of the symplectic normal connected sum $X_1(1)$, then we obtain a minimal symplectic $4$-manifold $N$\/ whose fundamental group is $\Z$.  $\pi_1(N)$ is generated by the image of the generator $c$\/ from $\pi_1(Y_1(1,1)\setminus\nu\Sigma_2)$.  In fact, $N=\tilde{X}_1(0,1)$ in the notation of Section~\ref{sec:exotic 1-2}.  Let $T'=a_2' \times c' \subset N$.  It is clear that the inclusion induced homomorphism $\pi_1(T')\rightarrow \pi_1(N)$ is surjective.  A meridian of $T'$ is given by the commutator $[d^{-1},b_2^{-1}]$, which can easily be shown to be trivial in $\pi_1(N\setminus\nu T')$ using exactly the same argument as in the proof of Lemma~\ref{lem:1-2}.  Thus we have constructed a pair $(N,T')$ corresponding to $(\chi,c)=(1,7)$.  

The choices for the other $(\chi,c)$ pairs can be verified in similar ways.  Finally, by perturbing the symplectic form on $N$, we can turn the Lagrangian torus $T'$ into a symplectic submanifold of $N$.  This concludes the proof of Theorem~\ref{thm:wedge}.

\section*{Acknowledgments}  
The authors thank Selman Akbulut, John B. Etnyre, Peter S. Ozsv\'{a}th and Zolt\'{a}n Szab\'{o} for their kind encouragements.  The authors also thank R. \.{I}nan\c{c} Baykur, Ronald Fintushel, Rafael Torres and the referee for very helpful comments.  
A. Akhmedov was partially supported by NSF grant FRG-0244663.  
B. D. Park was partially supported by CFI, NSERC and OIT grants.

\end{document}